\documentclass[12pt]{article}

\def\b{{\bf b}}
\def\bl{\begin{eqnarray}}
\def\el{\end{eqnarray}}
\def\bll{\begin{eqnarray*}}
\def\ell{\end{eqnarray*}}
\usepackage{amssymb}
\usepackage{amsmath,amsthm}
\usepackage{amsfonts}
\usepackage{fleqn}
\usepackage{showidx}
\usepackage{titlesec}
\usepackage{makeidx}
\makeindex
\usepackage[backref=true,dvipdfm,draft=false,colorlinks=true,
bookmarksopen=true,bookmarks=true,bookmarksnumbered=false,
linkcolor=blue,citecolor=blue,urlcolor=blue,pdfauthor=YunBaoHuang]{hyperref}

\title{The powers of smooth words over arbitrary 2-letter alphabets}

\small \author{Yun Bao Huang \\
 Department of Mathematics\\
 Hangzhou Normal University\\
 Xiasha Economic Development Area\\
 Hangzhou, Zhejiang 310036, China\\
 huangyunbao@sina.com\\
 huangyunbao@gmail.com}
\date{2011.08.07}

\begin{document}
\numberwithin{equation}{section}
\makeatletter
\titlelabel{\thetitle.\;}
\newcommand{\Extend}[5]{\ext@arrow 0099{\arrowfill@#1#2#3}{#4}{#5}}
\newcommand{\prf}[1]{\noindent\bf{Proof.}\;\;\rm{#1}{.\;$\Box$}}
\newtheorem{defn}{Definition}
\newtheorem{thm}[defn]{Theorem}
\newtheorem{lem}[defn]{Lemma}
\newtheorem{prop}[defn]{Proposition}
\newtheorem{cor}[defn]{Corollary}
\newtheorem{conj}[defn]{Conjecture}
\newtheorem{exmp}[defn]{Example}
\newtheorem{rem}[defn]{Remark}
\makeatother \maketitle
\newpage
\begin{quote}
{\small {\bf Abstract.} Carpi (1993) and Lepist\"{o} (1994) proved
independently that smooth words are cube-free for the alphabet
$\{1,2\}$, but nothing is known on whether for the other 2-letter
alphabets, smooth words are $k$-power-free for some suitable
positive integer $k$. In this paper, we first establish the
derivative formula of the concatenation of two smooth words and
power derivative formula of smooth words for arbitrary 2-letter
alphabets. Then, by making use of power derivative formula, for
arbitrary 2-letter alphabet $\{a,\,b\}$ with $a\text{ and }b$ being
positive integers and $a<b$, we prove that smooth words are
$(b+1)$-power-free except for $a=1 \text{ and }b=3$, and smooth
words are quintic-free and there are infinitely many smooth
biquadrates for the alphabet $\{1,\,3\}$. Moreover, we give the
number $\gamma_{a,b}(n)$ of smooth words of the form $w^{n}$ for
$a$ and $b$ having same parity.\\
{\bf Keywords:} Closure; Differentiable; Closurely differentiable;
Derivative; Smooth words; Power derivative formula; Power-free
index.}
\end{quote}
\newpage
\section{Introduction}
The curious Kolakoski sequence $K$~\cite{Kolakoski} is the infinite
sequence over the alphabet $\Sigma=\{1,2\}$, which starts with 1 and
equals the sequence defined by its run lengths:

$K=
\underbrace{1}_{1}\underbrace{22}_{2}\underbrace{11}_{2}\underbrace{2}_{1}
 \underbrace{1}_{1}\underbrace{22}_{2}\underbrace{1}_{1} \underbrace{22}_{2}\underbrace{11}_{2}
 \underbrace{2}_{1}\underbrace{11}_{2}\underbrace{22}_{2}\underbrace{\cdots}_{\cdots\;\;=\;K}$\\
Here, a run is a maximal subsequence of consecutive identical
symbols. The Kolakoski sequence $K$ has received a remarkable
attention by showing some intriguing combinatorical properties (see
Kimberling~\cite{Kimberling}, Dekking~\cite{Dekking2,Dekking3}), one
of them is known as Keane's problem~\cite{Keane}:

Is the frequency of 1(or 2) in the sequence $K$  equal to
$\frac{1}{2}$ ?

Chv\'{a}tal \cite{Chv} made use of a clever approach to obtain that
the letter frequencies are limited to $0.5\pm 0.000838$. Moreover,
many other studies have been done on $K$ in Brlek et al.
\cite{Br1,Br2,Br3,Br4,Br5}, Huang~\cite{Huang1,Huang2,Huang4}, Huang
and Weakley~\cite{Huang3}, Steacy~\cite{Steacy},
Steinsky~\cite{Steinsky} and Weakley~\cite{Weakley}. Especially,
Baake and Sing~\cite{Ba1} and Sing~\cite{Sing2, Sing3} established a
connection between generalized Kolakoski sequences and model sets.
Sing~\cite{Sing4} gave a recent survey on Kolakoski sequences and
related problems.

   P\v{a}un~\cite{P} conjectured that the Kolakoski sequence $K$
contains only squares of bounded length and that it is cube-free.
Carpi~\cite{Carpi1} and  Lepist\"{o}~\cite{Lepist} independently
solved the conjecture respectively.
   Carpi~\cite{Carpi2} further showed that for any positive integer n, only
finitely many words can occur twice, at distance n, in a
$C^{\infty}$-word.

A naturally arising question is whether or not for the other
2-letter alphabets, smooth words are $k$-power-free for some
suitable positive integer $k$. This paper is a study of the property
of smooth words being $k$-power-free for arbitrary 2-letter
alphabets. However, an extension of the attractive result, which
smooth words are cube-free for the alphabet $\{1,2\}$, to arbitrary
2-letter alphabets leads to difficulties if we attempt to follow
Carpi's~\cite{Carpi2} or Lepist\"{o}'s~\cite{Lepist} method. But
surprisingly, the power derivative formula provides us the required
method to establish the power-free index of smooth words over
arbitrary 2-letter alphabets, which do not rely on machine
computation.

The paper is structured as follows. In Section \ref{s1}, we shall
first fix some notations and introduce some notions. Second in
Section \ref{s2}, we establish the derivative formula of the
concatenation of two smooth words (Theorem \ref{T1}) and power
derivative formula of smooth words (Theorem \ref{T3}) on any
2-letter alphabets.  In Section \ref{s4}, for arbitrary alphabet
$\{a,b\}$, we prove that all smooth words are
$\delta(a,b)$-power-free, where if $a=1\text{ and } b=3$, then
$\delta(a,b)=b+2$, otherwise, $\delta(a,b)=b+1$.  In Section
\ref{s5}, we give the number of smooth words of the form $u^{n}$ for
2-letter alphabets having same parity, where $n$ is a positive
integer. Finally, in Section \ref{s7}, we give some problems on
smooth words which are deserved to attention.
\section{Definitions and notation\label{s1}}
Let $\Sigma=\{a, b\}$ with $a<b$ and $a,\,b$ being positive
integers, $\Sigma^{*}$ denotes the free monoid over $\Sigma$ with
$\varepsilon$ as the empty word, and $\Sigma^{+}$ denotes
$\Sigma^{*}-\{\varepsilon\}$. A \textit{finite word} over $\Sigma$
is an element of $\Sigma^{*}$. If $w=w_{1}w_{2}\cdots$$w_{n}$,
$w_{i}\in \Sigma$ for $i=1, 2, \cdots, n$, then $n$ is called the
\textit{length} of the word $w$ and is denoted by $|w|$.

Let $\mathbb{N}$ be the set of positive integers. Then an
\textit{infinite word} on $\Sigma$ is a function
$w:\,\mathbb{N}\rightarrow \Sigma,$ that is $w=w_1w_2\cdots
w_n\cdots$. The set of all infinite words is denoted by
$\Sigma^{\omega}$. Given a word $w\in\Sigma^{*}$, a
$\mathit{factor}$ (or $\mathit{subword}$) $u$ of $w$ is a word $u\in
\Sigma^{*}$ such that $w=xuy$ for $ x,\, y \in\Sigma^{*}$.
   A $\mathit{run}$
(or $\mathit{block}$) is  a maximal factor of consecutive same
letters. Finally, $\mathbb{N}^{*}$ denotes the free monoid over
$\mathbb{N}$, and $\mathbb{N}^+=\mathbb{N}-\{\varepsilon\}$.

   The \textcolor[rgb]{0.00,0.00,1.00}{$\mathit{reversal}$} (or $\mathit{mirror\, image}$)\label{def4} of $u=u_{1}u_{2}\cdots$$u_{n}\in\Sigma^{*}$ is
the word $\tilde{u}=u_{n}u_{n-1}\cdots$$u_{2}\\u_{1}$. The
\textcolor[rgb]{0.00,0.00,1.00}{$\mathit{complement}$} (or
$\mathit{permutation}$) of $u=u_{1}u_{2}\cdots$$u_{n}\in\Sigma^{*}$
is the word $\bar{u}=\bar{u}_{1}\bar{u}_{2}\cdots$$\bar{u}_{n}$,
where $\bar{a}=b, \bar{b}=a$.

  Every word $w\in\Sigma^+$ can be uniquely written as a
product of factors as follows:

   $w=\alpha^{i_{1}}\bar{\alpha}^{i_{2}}\alpha^{i_{3}}\bar{\alpha}^{i_{4}}\cdots$, where $i_{j}>0$, $j\in \mathbb{N}$.\\
The operator giving the size of the blocks appearing in the coding
is a function:

  $\Delta:\Sigma^+\rightarrow N^+$\label{def5}, defined by

  $\Delta(w)=i_{1}i_{2}i_{3}\cdots=\prod_{k\geq 1}i_{k}$\\
which is easily extended to infinite words.

  For any $w\in \Sigma^{*}$ (or $\Sigma^{\omega})$,
  $\mathit{first(w)}$
denotes the first letter of the word $w$. For each $w \in
\Sigma^{*}$, $\mathit{last(w)}$ denotes the last letter of the word
$w$. It is clear that the operator $\Delta$ satisfies the property:
   $\Delta(uv)=\Delta(u)\Delta(v)$ if and only if  $\mathit{last(u)\neq
   first(v)}$.

Pseudo-inverse functions:

  $\Delta_{a}^{-1}, \Delta_{b}^{-1}:
N^+\rightarrow\Sigma^+,\quad u=u_{1}u_{2}u_{3}u_{4}\cdots$\\
 can be defined by

  $\Delta_{\alpha}^{-1}(u)=\alpha^{u_{1}}\bar{\alpha}^{u_{2}}\alpha^{u_{3}}\bar{\alpha}^{u_{4}}\cdots$\\
where $\alpha\in \Sigma$.

   For $w\in \Sigma^{*}$,  $r(w)$ denotes the number of runs of $w$,
$fr(w)$ and $lr(w)$ denote the first run and last run of $w$
respectively, and $lfr(w)$ and $llr(w)$ denote the length of the
first run and last run of $w$ respectively. For example, if
$w=a^{2}b^{2b}a^{4a}b^{3}$, then $fr(w)=a^{2}$, $lr(w)=b^{3}$,
$lfr(w)=2$ and $llr(w)=3$. Next we introduce the notion of the
closure of a word $w\in \Sigma^{*}$.
\begin{defn}\label{def1} \rm
Let $w\in \Sigma^{*}$ and
\begin{equation}
w=\alpha^{t_{1}}\bar{\alpha}^{t_{2}}\ldots \beta^{t_{k}},
\textrm{where} \,\alpha,\,\beta\in \Sigma, \,1\leq t_{i}\leq b \,\, \textrm{for} \,1\leq i\leq k; \label{hk1}\\
\end{equation}
$\quad\quad\;\; \hat{w}=\left\{\begin{array}{ll}
w ,& lfr(w)\leq a \;\textrm{and}\;llr(w)\leq a\\
\alpha^{b-t_{1}}w,& lfr(w)> a \;\textrm{and}\;llr(w)\leq a\\
w\beta^{b-t_{k}}, & lfr(w)\leq a \;\textrm{and}\;llr(w)> a\\
\alpha^{b-t_{1}}w\beta^{b-t_{k}}, &  lfr(w)> a \;\textrm{and}\;llr(w)> a\\
\end{array}\right..$\\
Then $\hat{w}$ is said to be the closure of $w$.
\end{defn}
For example, let $w=3311133313133311133,\, u=3313133311$, then $u$
is a factor of $w$, and $\hat{w}=333111333131333111333,\,
\hat{u}=333131333111$. Thus $\hat{u}$ is a factor of $\hat{w}$,
which also holds in general (see Lemma \ref{Lem2} (1)).

Now we generalize the definition of differentiable words, which
Dekking introduced in~\cite{Dekking2},  to an arbitrary 2-letter
alphabet $\{a,\,b\}$ from the alphabet $\{1,\,2\}$ .

\begin{defn}\label{def2}\rm
Let $w\in \Sigma^{*}$ be of the form  (\ref{hk1}). If the length of
every run of $w$ takes only $a\,\textrm{or}\, b$ except for the
length of first and last runs, then we call that $w$ is
differentiable, and its derivative, denoted by $D(w)$, is the word
whose $j$th symbol equals the length of the $j$th run of $w$,
discarding the first and/or the last run if its length is less than
$b$.
\end{defn}
If $\hat{w}$ is differentiable, then we call that $w$ is closurely
differentiable. If a finite word $w$ is arbitrarily often closurely
differentiable, then we call $w$ a $C^{\infty}_{a,b}$-word or a
smooth word over the alphabet $\{a,b\}$, and the set of all smooth
words over the alphabet $\{a,b\}$ is denoted  by $C^{\infty}_{a,b}$
or $C^{\infty}$.

Let $\rho(w)=D(\hat{w})$\label{def6}, then it is clear that $w$ is a
smooth word if and only if there is a positive integer $k$ such that
$\rho^{k}(w)=\varepsilon$.

Moreover, it is clear that D is an operator from $\Sigma^{*}$ to
$\Sigma^{*}$, $r(w)\leq |D(w)|+2$\label{huang} and
\begin{eqnarray}\label{hk2}
D(\hat{w})=\left\{\begin{array}{ll}
bD(w),& b>lfr(w)> a \;\textrm{and}\;llr(w)\leq a\\
D(w)b, &  b>llr(w)> a \,\;\textrm{and}\;lfr(w)\leq a\\
bD(w)b, &  b>lfr(w)> a \;\textrm{and}\;b>llr(w)> a\\
D(w) ,& \textrm{otherwise}
\end{array}\right..
\end{eqnarray}
By  (\ref{hk2}), it is obvious that if $w$ is closurely
differentiable, then it must be differentiable.

   Obviously, if $w$ is a smooth word and $|w|>0$, then  $|D(w)|<|w|$.
Moreover, $D$ and $\Delta$ can be both iterated.
\begin{defn}\label{def3}\rm A word $w\in \Sigma^{\omega}$ is a \it $C^{\omega}_{a,b}$-word \rm or \it
smooth infinite word \rm if $\Delta^{k}(w)\in \Sigma^{\omega}$ for
all $k$ $\in N$. The class of smooth infinite words is denoted by
$C_{a,b}^{\omega}\text{ or }C^{\omega}$.
\end{defn}

It is easy to see that finite factors of $C^{\omega}_{a,b}$-words
are all $C^{\infty}_{a,b}$-words. Thus finite smooth
words~\cite{Br1} are $C^{\infty}_{a,b}$-words.

It is hardly effortless to check that $D$ commute with the mirror
image ($\tilde{\,\,}$) and are stable for the
permutation($\bar{\,\,}$) over arbitrary 2-letter alphabet
$\{a,b\}$. Thus Lemmas 1-2 in~\cite{Huang5} still holds over
arbitrary 2-letter alphabets $\{a,b\}$.
\begin{lem}\rm
\it\; For each $u\in\Sigma^{*},\; D(\tilde{u})=\widetilde{D(u)},\;
D(\bar{u})=D(u)$.
\end{lem}
\section{Derivative Formula\label{s2}}
The following Lemmas \ref{Lem2}-\ref{Lem4}  are important in what
follows.
\begin{lem}[~\cite{Huang4} Lemma 5]\label{Lem2}
Let $w$ be a differentiable word  and $u$ is a factor of $w$. Then

\rm (1) \it $\hat{u}$ and $w$ are both factors of $\hat{w}$;

\rm (2) \it $\hat{\tilde{w}}=\tilde{\hat{w}},
\hat{\bar{w}}=\bar{\hat{w}}$;

\rm (3) \it $D(u)$ is a factor of $D(w)$, $D(w)$ is a factor of
$\Delta(w)$;

\rm (4) \it If $w$ is closurely differentiable, then $\rho(u)$ and
$D(w)$ are both factors of $\rho(w)$, and $\rho(\bar{w})=\rho(w)$,
$\rho(\tilde{w})=\widetilde{\rho(w)}$.\;$\Box$
\end{lem}

Lemma \ref{Lem2} (4) means that $C^\infty$ is closed under the
operators mirror image and complementation:

  $w \in C^{\infty}\Longleftrightarrow\bar{w}, \tilde{w}\in
C^{\infty}$.

\begin{lem}[\cite{Huang4} Lemma 7]\label{Lem4}\rm
(1) \it Let $w=w_{1}w_{2}\cdots w_{n}$ be a
 smooth word. Then any factor of $w$ is also a
 smooth word;\rm

 (2) \it Any smooth word $w=w_{1}w_{2}\cdots w_{n}$ has both a left smooth extension and a right
smooth extension;\rm

 (3) \it If \, $w\in \Sigma^{*}$ and $\Delta(w)\in C^{\infty}$, then $w\in
C^{\infty}$.\;$\Box$
\end{lem}

   Huang~\cite{Huang5} obtained the derivative
formula of powers of a $C^\infty$-word for the alphabet $\{1,2\}$.
Now we establish the derivative formula of the concatenation of
two smooth words for arbitrary 2-letter alphabet.
\begin{thm}[Derivative Formula]\label{T1}\;
For any $x \in D_{\Sigma}$, if $uxv \in C^{\infty}$ then there
exists an element $w \in D_{\Sigma}$ such that $D(uxv)=D(u)wD(v)$,
where $D_{\Sigma}=D_{a,b}$ and
\begin{eqnarray}
D_{1,2}&=&\{\varepsilon,1,2,12,21,11,22,
112,211,121,122,221,212,1121,1211,1212,
                                   \nonumber\\
 & & 2121,2112, 1221,1122,2211,11211\};\label{maker4}\nonumber
\end{eqnarray}
\vspace{-1cm}
\begin{eqnarray}
D_{1,3}  =  \{\varepsilon,1,3,13,31,11,33,
113,311,131,313,111,3111,1113,1311,1131\}; \label{maker5}\nonumber
\end{eqnarray}
\begin{eqnarray}
D_{1,4}  =  \{\varepsilon,1,4,14,41,11,44,
111,411,114,141,414,1111,4111,1114\}; \label{maker6}\nonumber
\end{eqnarray}
\begin{eqnarray}
D_{1,b}  =  \{\varepsilon,1,b,1b,b1,11,bb,
11b,b11,111,1111\},\;\text{where}\; b\geq 5;
\label{maker7}\nonumber
\end{eqnarray}
\begin{eqnarray}
D_{2,b} =  \{\varepsilon,2,b,2b,b2,22,bb,222\};
\label{maker8}\nonumber
\end{eqnarray}
\begin{eqnarray}
D_{a,b}  =  \{\varepsilon,a,b,aa,bb,ab,ba\},\;
\textrm{where}\;a\geq 3. \label{maker9}\nonumber
\end{eqnarray}
\end{thm}
\noindent{\bf Proof.} We divide the proof  of Theorem \ref{T1} into
the six cases according to the value of $D_{\Sigma}$. In what
follows, $\alpha\in \Sigma$. We easily see \vspace{-.33cm}
\begin{equation}
D(uxv)=D(u)\Delta(yxz)D(v)\text{ where }y,z\in\{\varepsilon,
a,a^2,\ldots,a^{b-1}, b,b^2,\ldots,b^{b-1} \}\;\label{eq1}
\end{equation}
{\bf Case 1.} $\Sigma=\{1,\,2\}$. If $|x|\leq 2$ then from
(\ref{eq1}) we have $|yxz|\leq 4$, which means that
$|\Delta(yxz)|\leq 3$, that is, Theorem~\ref{T1} holds; If $|x|=3$
then
$x=\alpha\bar{\alpha}\alpha,\alpha\alpha\bar{\alpha},\alpha\alpha\bar{\alpha}\alpha\alpha$.
It is easy to verify that $\Delta(yxz)\in D_{1,2}$; If $|x|=4$ then
$x=\alpha\bar{\alpha}\alpha\bar{\alpha},\alpha\alpha\bar{\alpha}\alpha,\alpha\alpha\bar{\alpha}\bar{\alpha},\text{
or } \\\alpha\bar{\alpha}\bar{\alpha}\alpha$, it is easy to see
$|\Delta(yxz)|\leq 4$ and  $\Delta(yxz)\neq 2212, 2122$ for the
first three cases; and $|\Delta(yxz)|\leq 4$ and  $\Delta(yxz)\neq
2212, 2122$ or $\Delta(yxz)=11211$ for the last case, which suggests
that Theorem~\ref{T1} also holds; If $x=11211$ then
$\Delta(yxz)=212,1212,2121$, that is, Theorem~\ref{T1}
holds.\\
{\bf Case 2.} $\Sigma=\{1,\,3\}$. If $|x|\leq 2$ then from
(\ref{eq1}) we have $|yxz|\leq 6$ and $yxz$ contains at most three
runs. Thus  $|\Delta(yxz)|\leq 3$ and if $\Delta(yxz)$ has three
runs then $\Delta(yxz)$ contains at most one run of length 3, which
means $\Delta(yxz)\in \Sigma_{1,3}$; If $|x|=3$ then it is easy to
get $|\Delta(yxz)|\leq 3$ and  $yxz$ contains at most one run of
length 3 except for $yxz=3^313^3$, which implies that
$\Delta(yxz)\in \Sigma_{1,3}$; If $|x|=4$, then by a direct check we
can obtain $\Delta(yxz)\in \Sigma_{1,3}$.

The remaining four cases can also be verified similarly. $\Box$

Now we generalize Theorem \ref{T1} to the following more general
form.
\begin{thm}\label{T2}
Let $u_{i}\in C^{\infty}$ for $i=1,2,3,\cdots,n\;(\geq 2)$. Then for
any positive integer $k$, there exist $w_{j}\in D_{\Sigma}$ for
$j=1,2,\cdots, n-1$ such that $D^{k}(u_{1}u_{2}\cdots
u_{n})=D^{k}(u_{1})w_{1}D^{k}(u_{2})w_{2}\cdots
D^{k}(u_{n-1})w_{n-1} D^{k}(u_{n})$.
\end{thm}
\noindent{\bf Proof.} First of all, we prove that the assertion
holds for n=2. For this, we only need to proceed by induction on
$k$. If $k=1$, in view of $u_{1}u_{2}\in C^{\infty}$ and taking
$x=\varepsilon$ in Theorem \ref{T1}, one sees that Theorem \ref{T2}
holds for $k=1$.

Now we suppose that Theorem \ref{T2} holds for all $k\leq m\;(\geq
1)$,
   i.e. there exists a $x_{1}\in D_{\Sigma}$ such that
   $D^{m}(u_{1}u_{2})=D^{m}(u_{1})x_{1}D^{m}(u_{2})$. Thus
   $D^{m+1}(u_{1}u_{2})=D(D^{m}(u_{1}u_{2}))=D(D^{m}(u_{1})x_{1}D^{m}(u_{2}))$, in view of
   Theorem \ref{T1}, one sees that there is a $w_{1}\in D_{\Sigma}$ such
   that $D(D^{m}(u_{1})x_{1}D^{m}(u_{2}))=D^{m+1}(u)w_{1}D^{m+1}(v)$, which implies
   that Theorem \ref{T2} holds for $n=2$.

From the above it immediately follows that Theorem \ref{T2} holds
for any positive integer $n\geq 2$.
 $\Box$

Note that if $u$ has at least two runs, then $D(u^{2})=D(u)xD(u)$,
where $x$ is uniquely determined by the $\mathit{fr(u)}$ and
$\mathit{lr(u)}$. Thus from Theorem \ref{T2}, we obtain the
following important result.
\begin{thm}[Power Derivative Formula]\label{T3}\;
Let $u\in C^{\infty}_{a,b}$, $n\;(\geq 2)$ and $k$ be  positive
integers. If $D^{k-1}(u)$ has at least two runs, then for any
positive integer $1\leq j\leq k$, there exists a word $w\in D_{a,b}$
such that $D^{j}(u^{n})=(D^{j}(u)w)^{n-1}D^{j}(u)$.
\end{thm}
\section{The power-free index of smooth words\label{s4}}
Let $n$ be a positive integer, if  $u^{n}\notin C^{\infty}_{a,b}$
for any nonempty $C^{\infty}_{a,b}$-word $u$, then we call that
$C^{\infty}_{a,b}$-words are $n$-power-free. And the minimal
positive integer $n$ such that $C^{\infty}_{a,b}$-words are
$n$-power-free is said to be the power-free index of
$C^{\infty}_{a,b}$-words.

Carpi~\cite{Carpi1} and Lepist\"{o}~\cite{Lepist} independently
proved the following interesting result.
\begin{prop}[Carpi and Lepist\"{o}]\label{p2}
The power-free index of $C^\infty_{1,2}$-words is equal to 3, that
is, smooth words are cube-free for the alphabet $\{1,2\}$.
\end{prop}
In this section, we establish the power-free index of smooth words
for arbitrary 2-letter alphabet $\{a,\,b\}$ that do not rely on
machine computation.
\begin{thm}\label{T4}
Let
\begin{eqnarray}
h(a,b)&=&\left\{\begin{array}{ll}
b+2 ,& a=1\,\textrm{ and }\,b=3\\
\frac{b+4}{2} ,& 2\mid b\\
\frac{b+5}{2},& 2\nmid b,\,\textrm{ and }\,b\neq 3,\;a=1\\
\frac{b+3}{2}, & 2\nmid b, \,\textrm{ and }\,a\geq 2
\end{array}\right.,\label{e5.1}\\
\delta(a,b)&=&\left\{\begin{array}{ll}
b+2 ,& a=1\,\textrm{ and }\,b=3\\
b+1 ,& \textrm{otherwise}
\end{array}\right.,\label{e5.2}
\end{eqnarray}
then\rm

 (1) \it Smooth words of length more than 1 are $h(a,b)$-power-free
for arbitrary $2$-letter alphabet $\{a,\,b\}$;\rm

 (2)  \it The power-free index of smooth words
is $\delta(a,b)$ for the alphabet $\{a,\,b\}$.
\end{thm}
\noindent{\bf Proof.} It is clear that $h(a,b)\geq 3$. Note that
$\widetilde{D_\Sigma}=D_\Sigma$. We easily see that the following
assertion holds.
\begin{eqnarray}
&&(\alpha^{i}x)^{h(a,b)-1}\alpha^{i}\in
C^{\infty}_{a,b}\Longleftrightarrow
(\alpha^{i}\tilde{x})^{h(a,b)-1}\alpha^{i}\in C^{\infty}_{a,b}.
\end{eqnarray}
{\bf (1)} By Proposition \ref{p2}, it suffices to check the
assertion (1) for the other cases except for $a=1\text{ and }b=2$.\\
{\bf Case 1.} $\Sigma=\{1,3\}$. Assume on the contrary that there
exists a smooth word $u$ of length more than or equal to 1  such
that $u^{5}\in C^{\infty}_{1,3}$. Then $u$ has at least two runs.
Let $k$ be the maximal integer such that $D^{k-1}(u)$ has at least
two runs. Thus by Power Derivative Formula, we have
\begin{equation}
(D^{j}(u)x)^{4}D^{j}(u)=D^{j}(u^{5})\in C^{\infty}_{1,3}\;
\textrm{for}\,1\leq j\leq k,\; \textrm{where}\; x\in
D_{1,3},\label{eq:maker171}
\end{equation}
and $D^{k}(u)$ has at most one run, which implies that

\quad$D^{k}(u)=\alpha^{i}$, where $i=0,1,2,3.$\\
{\bf Case 1.1.} $D^{k}(u)=\varepsilon$. Then since $D^{k-1}(u)$ has
at least two runs, we have
\begin{equation}
D^{k-1}(u)=\alpha^{s}\bar{\alpha}^{t},\;\textrm{where}\;s,t=1,2,\label{eq:maker172}
\end{equation}
and by  (\ref{eq:maker171}) ($j=k$), we can get
\begin{equation}
D^{k}(u^{5})=x^{4},\;\textrm{where}\;x\in
D_{1,3}.\label{eq:maker173}
\end{equation}
From  (\ref{eq:maker173}), a direct verification leads to

\quad $D^{k}(u^{5})=\varepsilon$.\\
Note that $s,\,t=1,2$, by (\ref{eq:maker172}) and
(\ref{eq:maker171}) ($j=k-1$), we see that $D^{k-1}(u^{5})$ has at
least ten runs, which means that $D^{k}(u^{5})$ has at least eight
runs, a contradiction to
$D^{k}(u^{5})=\varepsilon$.\\
{\bf Case 1.2.} $D^{k}(u)=\alpha^{i}$, $i$=1, 2, 3. Then by
(\ref{eq:maker171}), we obtain
$(\alpha^{i}x)^{4}\alpha^{i}=D^{k}(u^{5})\in C_{1,3}^{\infty}$,
where $x\in D_{1,3}$. If $i=1$ then $(\alpha x)^{4}\alpha\in
C^\infty$, which implies that $x\neq \alpha^j,\bar{\alpha}^t$, where
$0\leq j,\,t\leq 3$. Thus $x=\alpha y\bar{\alpha}\text{ or
}\bar{\alpha} y\alpha$, that is, $(\alpha^2
y\bar{\alpha})^4\alpha\text{ or }\alpha (\bar{\alpha}y\alpha^2)^4\in
C^\infty$. Then it follows that $y=\alpha z\text{ or }z\alpha$. So
$(\alpha^3 z\bar{\alpha})^4\alpha\text{ or }\alpha
(\bar{\alpha}z\alpha^3)^4\in C^\infty$, which suggests that
$z=\varepsilon\text{ or }\bar{\alpha}$ by $x\in D_{1,3}$, a
contradiction. If $i=2,3$ then by a similar argument we see that
there is no $x\in D_{1,3}$ such that $(\alpha^{i}x)^{4}\alpha^{i}\in
C_{1,3}^{\infty}$.\\
{\bf Case 2.} $\Sigma=\{1,4\}$. Suppose to the contrary that there
exists a smooth word  $u$ of length larger than 1 such that
$u^{4}\in C^{\infty}_{1,4}$. Then $u$ has at least two runs. Let $k$
be the maximal integer such that $D^{k-1}(u)$ has at least two runs.
Then by Power Derivative Formula, we obtain
\begin{equation}
(D^{j}(u)x)^{3}D^{j}(u)=D^{j}(u^{4})\in C^{\infty}_{1,4}\;
\textrm{for}\,1\leq j\leq k,\; \textrm{where}\; x\in
D_{1,4},\label{eq:maker68}
\end{equation}
and $D^{k}(u)$ has at most one run, which implies that

\quad$D^{k}(u)=\alpha^{i}$, where $i=0,1,2,3,4$.\\
{\bf Case 2.1.} $D^{k}(u)=\varepsilon$. Then since $D^{k-1}(u)$ has
at least two runs, we have
\begin{equation}
D^{k-1}(u)=\alpha^{s}\bar{\alpha}^{t},\;\textrm{where}\;s,t=1,2,3,\label{eq:maker69}
\end{equation}
and by  (\ref{eq:maker68}) we have

\quad $D^{k}(u^{4})=x^{3},\;\textrm{where}\;x\in D_{1,4}$.\\
A direct verification leads to
\begin{equation}
 D^{k}(u^{4})=\varepsilon,1^{3},4^{3},(14)^{3},(41)^{3}.\label{eq:maker70}
\end{equation}
From  (\ref{eq:maker68}) it follows that
\begin{equation}
(D^{k-1}(u)x)^{3}D^{k-1}(u)=D^{k-1}(u^{4}).\label{eq:maker71}
\end{equation}
Since (\ref{eq:maker69}) and (\ref{eq:maker70}) satisfy
(\ref{eq:maker71}), we see that $D^{k-1}(u^{4})$ has at least eight
runs, which means that $ D^{k}(u^{4})=(14)^{3},(41)^{3}$ by
(\ref{eq:maker70}).  Thus from (\ref{eq:maker71}) we obtain
\begin{eqnarray}
(\alpha^{s}\bar{\alpha}^{t}x)^{3}\alpha^{s}\bar{\alpha}^{t}
&=&\alpha^{i}\bar{\alpha}\alpha^{4}\bar{\alpha}\alpha^{4}\bar{\alpha}\alpha^{4}\bar{\alpha}^{j},\;1\leq i\leq 3,\,0\leq j\leq 3;\label{eq:maker72}\\
&=&\alpha^{i}\bar{\alpha}^{4}\alpha\bar{\alpha}^{4}\alpha\bar{\alpha}^{4}\alpha\bar{\alpha}^{j},\;0\leq
i\leq 3,\,1\leq j\leq 3.\label{eq:maker73}
\end{eqnarray}
{\bf Case 2.1.1.} If  (\ref{eq:maker72}) holds, then from
(\ref{eq:maker69}) we obtain
\begin{equation}
D^{k-1}(u)=\alpha^{i}\bar{\alpha}$  and $x=\alpha^{4-i}\in
D_{1,4},\, i=1,2,3\label{eq:maker74}
\end{equation}
{\bf Case 2.1.1.1.} i=1. Then by $x=\alpha^{3}\in D_{1,4}$, we have
$\alpha=1$. Thus from (\ref{eq:maker74}) and (\ref{eq:maker72}) we
obtain
\begin{eqnarray}
D^{k-1}(u)&=&14;\label{eq:maker75}\\
D^{k-1}(u^{4})&=&141^{4}41^{4}41^{4}4.\label{eq:maker76}
\end{eqnarray}
From (\ref{eq:maker75}) and (\ref{eq:maker76}) it follows that
\begin{eqnarray}
D^{k-2}(u)&=&\beta^{k}\bar{\beta}\beta^{4}\bar{\beta}^{l},\, \,\,1\leq k\leq 3,\,0\leq l\leq 3;\label{eq:maker77}\\
D^{k-2}(u^{4})&=&\gamma^{s}\bar{\gamma}\gamma^{4}\bar{\gamma}\gamma\bar{\gamma}\gamma\bar{\gamma}^{4}\gamma
\bar{\gamma}\gamma\bar{\gamma}\gamma^{4}\bar{\gamma}\gamma\bar{\gamma}\gamma\bar{\gamma}^{4}\gamma^{t},
\,1\leq s\leq 3,\,0\leq t\leq 3.\label{eq:maker78}
\end{eqnarray}
By  (\ref{eq:maker68}), we have
\begin{equation}
D^{k-2}(u^{4})=(D^{k-2}(u)x)^{3}D^{k-2}(u),\;\text{where}\,x\in
D_{1,4}.\label{eq:maker79}
\end{equation}
Thus,  (\ref{eq:maker77}) means that the right side of
(\ref{eq:maker79}) has at least
four same runs of the form $\beta^{4}$, a contradiction to  (\ref{eq:maker78}).\\
{\bf Case 2.1.1.2.} i=2. Then by (\ref{eq:maker74}) and
(\ref{eq:maker72}), we have
\begin{eqnarray}
D^{k-1}(u)&=&\alpha^{2}\bar{\alpha};\label{eq:maker80}\\
D^{k-1}(u^{4})&=&\alpha^{2}\bar{\alpha}\alpha^{4}\bar{\alpha}\alpha^{4}\bar{\alpha}\alpha^{4}\bar{\alpha},\;\alpha=1,\,4.\label{eq:maker81}
\end{eqnarray}

If $\alpha=1$ then from (\ref{eq:maker80}) and (\ref{eq:maker81}),
we obtain
\begin{eqnarray}
D^{k-2}(u)&=&\beta^{k}\bar{\beta}\beta\bar{\beta}^{4}\beta^{l},\, \,1\leq k\leq 3,\,0\leq l\leq 3;\label{eq:maker82}\\
D^{k-2}(u^{4})&=&\gamma^{s}\bar{\gamma}\gamma\bar{\gamma}^{4}\gamma\bar{\gamma}\gamma\bar{\gamma}\gamma^{4}\bar{\gamma}
\gamma\bar{\gamma}\gamma\bar{\gamma}^{4}\gamma\bar{\gamma}\gamma\bar{\gamma}\gamma^{4}\bar{\gamma}^{t},
\,1\leq s\leq 3,\,0\leq t\leq 3.\label{eq:maker83}
\end{eqnarray}
An argument similar to Case 2.1.1.1 leads to a contradiction to
(\ref{eq:maker83}).

If $\alpha=4$ then from (\ref{eq:maker80}) and (\ref{eq:maker81}),
we obtain
\begin{eqnarray}
D^{k-2}(u)&=&\beta^{k}\bar{\beta}^{4}\beta^{4}\bar{\beta}\beta^{l}, \;0\leq k\leq 3,\,1\leq l\leq 3;\label{eq:maker84}\\
D^{k-2}(u^{4})&=&\gamma^{s}\bar{\gamma}^{4}\gamma^{4}\bar{\gamma}\gamma^{4}\bar{\gamma}^{4}\gamma^{4}\bar{\gamma}^{4}\gamma\bar{\gamma}^{4}
\gamma^{4}\bar{\gamma}^{4}\gamma^{4}\bar{\gamma}\gamma^{4}\bar{\gamma}^{4}\gamma^{4}\bar{\gamma}^{4}\gamma\bar{\gamma}^{t},\label{eq:maker85}
\end{eqnarray}
where $0\leq s\leq 3,\,1\leq t\leq 3$. Thus  (\ref{eq:maker84})
implies that the right side  of  (\ref{eq:maker79}) has at least
four same factors of the form
$\bar{\beta}^{4}\beta^{4}\bar{\beta}\beta$, a
contradiction to  (\ref{eq:maker85}).\\
{\bf Case 2.1.1.3.} i=3. Then from (\ref{eq:maker74}) and
(\ref{eq:maker72}), we obaitn
\begin{eqnarray}
D^{k-1}(u)&=&\alpha^{3}\bar{\alpha};\label{eq:maker86}\\
D^{k-1}(u^{4})&=&\alpha^{3}\bar{\alpha}\alpha^{4}\bar{\alpha}\alpha^{4}\bar{\alpha}\alpha^{4}\bar{\alpha},\;\alpha=1,\,4.\label{eq:maker87}
\end{eqnarray}

If $\alpha=1$ then by (\ref{eq:maker86}) and (\ref{eq:maker87}), we
have
\begin{eqnarray}
D^{k-2}(u)&=&\beta^{l}\bar{\beta}\beta\bar{\beta}\beta^{4}\bar{\beta}^{k},\,1\leq l\leq 3,\,0\leq k\leq 3;\label{eq:maker88}\\
D^{k-2}(u^{4})&=&\gamma^{s}\bar{\gamma}\gamma\bar{\gamma}\gamma^{4}\bar{\gamma}\gamma\bar{\gamma}\gamma\bar{\gamma}^{4}
\gamma\bar{\gamma}\gamma\bar{\gamma}\gamma^{4}\bar{\gamma}\gamma\bar{\gamma}\gamma\bar{\gamma}^{4}\gamma^{t},\,1\leq
s\leq 3,\,0\leq t\leq 3.\label{eq:maker89}
\end{eqnarray}
An argument similar to Case 2.1.1.1 leads to a contradiction to
(\ref{eq:maker89}).

If $\alpha=4$ then from (\ref{eq:maker86}) and (\ref{eq:maker87}),
we get
\begin{eqnarray}
D^{k-2}(u)&=&\beta^{k}\bar{\beta}^{4}\beta^{4}\bar{\beta}^{4}\beta\bar{\beta}^{l},\;0\leq k\leq 3,\,1\leq l\leq 3;\label{eq:maker90}\\
D^{k-2}(u^{4})&=&\gamma^{s}\bar{\gamma}^{4}\gamma^{4}\bar{\gamma}^{4}\gamma\bar{\gamma}^{4}\gamma^{4}\bar{\gamma}^{4}\gamma^{4}
\bar{\gamma}\gamma^{4}\bar{\gamma}^{4}\gamma^{4}\bar{\gamma}^{4}\gamma\bar{\gamma}^{4}\gamma^{4}\bar{\gamma}^{4}\gamma^{4}
\bar{\gamma}\gamma^{t},\label{eq:maker91}
\end{eqnarray}
where $0\leq s\leq 3,\,1\leq t\leq 3$. Therefore, (\ref{eq:maker90})
suggests that the right side  of  (\ref{eq:maker79}) has at least
four same factors of the form
$\bar{\beta}^{4}\beta^{4}\bar{\beta}^{4}\beta\bar{\beta}$, a
contradiction to  (\ref{eq:maker91}).\\
{\bf Case 2.1.2.} If  (\ref{eq:maker73}) holds, then from
(\ref{eq:maker69}), we obtain
\begin{eqnarray}
D^{k-1}(u)&=&\alpha\bar{\alpha}^{j}, \,1\leq j\leq 3;\nonumber\\
D^{k-1}(u^{4})&=&\alpha\bar{\alpha}^{4}\alpha\bar{\alpha}^{4}\alpha\bar{\alpha}^{4}\alpha\bar{\alpha}^{j},\;1\leq
j\leq 3,\;\alpha=1,\,4.\nonumber
\end{eqnarray}
Thus, by $\widetilde{D(w)}=D(\tilde{w})$, we have
\begin{eqnarray}
D^{k-1}(\tilde{u})&=&\bar{\alpha}^{j}\alpha, \,1\leq j\leq 3;\nonumber\\
D^{k-1}(\tilde{u}^{4})&=&\bar{\alpha}^{j}\alpha\bar{\alpha}^{4}\alpha\bar{\alpha}^{4}\alpha\bar{\alpha}^{4}\alpha,\;1\leq
j\leq 3,\nonumber
\end{eqnarray}
which means that $D^{k-1}(\tilde{u})$ and $D^{k-1}(\tilde{u}^{4})$
satisfy (\ref{eq:maker72}) and (\ref{eq:maker74}). So by Case
2.1.1,
we arrive at a contradiction.\\
{\bf Case 2.2.} $D^{k}(u)=\alpha$. Then by  (\ref{eq:maker68}), we
have

\quad $(\alpha x)^{3}\alpha=D^{k}(u^{4})\in C^{\infty}_{1,4}$,
where $x\in D_{1,4}$,\\
which implies that
\begin{eqnarray}
D^{k}(u)&=&\alpha,\label{eq:maker92}\\
D^{k}(u^{4})&=&\alpha^{4},\;\alpha=1,4;\label{eq:maker93}
\end{eqnarray}
or
\begin{eqnarray}
D^{k}(u)&=&4;\label{eq:maker94}\\
D^{k}(u^{4})&=&(41111)^{3}4;\label{eq:maker95}
\end{eqnarray}
or
\begin{eqnarray}
D^{k}(u)&=&1;\label{eq:maker96}\\
D^{k}(u^{4})&=&(11114)^{3}1;\label{eq:maker97}\\
            &=&1(41111)^{3}.\label{eq:maker98}
\end{eqnarray}
{\bf Case 2.2.1.} If (\ref{eq:maker92}) and (\ref{eq:maker93}) hold
then \\
{\bf Case 2.2.1.1.} $\alpha$=1. Then since (\ref{eq:maker92}) and
(\ref{eq:maker93}) satisfy  (\ref{eq:maker68}), we have
\begin{eqnarray}
D^{k-1}(u)&=&\beta^{i}\bar{\beta}\beta^{j},\;1\leq i,\,j \leq 3;\label{eq:maker99}\\
D^{k-1}(u^{4})&=&\beta^{i}\bar{\beta}\beta\bar{\beta}\beta\beta^{j},\;1\leq
i,\,j \leq 3,\;\beta=1,4.\label{eq:maker100}
\end{eqnarray}
Thus (\ref{eq:maker71}) and (\ref{eq:maker99}) mean that
$D^{k-1}(u^{4})$ has at least eight runs, contradicts
(\ref{eq:maker100}).\\
{\bf Case 2.2.1.2.} $\alpha$=4. Then similarly, we have
\begin{eqnarray}
D^{k-1}(u)&=&\beta^{i}\bar{\beta}^{4}\beta^{j},\;0\leq i,\,j \leq 3;\label{eq:maker101}\\
D^{k-1}(u^{4})&=&\beta^{i}\bar{\beta}^{4}\beta^{4}\bar{\beta}^{4}\beta^{4}\beta^{j},\;0\leq
i,\,j \leq 3.\label{eq:maker102}
\end{eqnarray}
Thus (\ref{eq:maker71}) and (\ref{eq:maker101}) mean that
$\bar{\beta}^{4}$ occurs at least four times in $D^{k-1}(u^{4})$,
contradicts  (\ref{eq:maker102}).\\
{\bf Case 2.2.2.} If (\ref{eq:maker94}) and (\ref{eq:maker95}) hold,
then we obtain
\begin{eqnarray}
D^{k-1}(u)&=&\alpha^{i}\bar{\alpha}^{4}\alpha^{j},\;0\leq i,\,j \leq 3\label{eq:maker103}\\
D^{k-1}(u^{4})&=&\alpha^{i}\bar{\alpha}^{4}\alpha\bar{\alpha}\alpha\bar{\alpha}\alpha^{4}
\bar{\alpha}\alpha\bar{\alpha}\alpha\bar{\alpha}^{4}\alpha\bar{\alpha}\alpha\bar{\alpha}\alpha^{4}\bar{\alpha}^{j},\;0\leq
i,\,j \leq 3\label{eq:maker104}
\end{eqnarray}
Hence, an argument similar to Case 2.2.1.2 gives rise to a contradiction.\\
{\bf Case 2.2.3.} If (\ref{eq:maker96}) and (\ref{eq:maker97}) hold
then (\ref{eq:maker99}) holds, and
\begin{equation}
D^{k-1}(u^{4})=\alpha^{i}\bar{\alpha}\alpha\bar{\alpha}\alpha\bar{\alpha}^{4}\alpha
\bar{\alpha}\alpha\bar{\alpha}\alpha^{4}\bar{\alpha}\alpha\bar{\alpha}\alpha\bar{\alpha}^{4}\alpha\bar{\alpha}^{j},\;1\leq
i,\,j \leq 3.\label{eq:maker105}
\end{equation}
Thus, since (\ref{eq:maker99}) and (\ref{eq:maker105}) satisfy
(\ref{eq:maker71}), which means that $\alpha=\beta$,  and the last
runs of (\ref{eq:maker99}) and (\ref{eq:maker105}) are the same,
which leads
to a contradiction.\\
{\bf Case 2.2.4.} If (\ref{eq:maker96}) and (\ref{eq:maker98}) hold
then (\ref{eq:maker99}) holds, and
\begin{equation}
D^{k-1}(u^{4})=\alpha^{i}\bar{\alpha}\alpha^{4}\bar{\alpha}\alpha\bar{\alpha}\alpha\bar{\alpha}^{4}\alpha
\bar{\alpha}\alpha\bar{\alpha}\alpha^{4}\bar{\alpha}\alpha\bar{\alpha}\alpha\bar{\alpha}^{j},\;1\leq
i,\,j \leq 3.\label{eq:maker106}
\end{equation}
Thus, (\ref{eq:maker99}) and (\ref{eq:maker106}) satisfy
(\ref{eq:maker71}), which suggests that $\alpha=\beta$,  and the
last runs of (\ref{eq:maker99}) and
(\ref{eq:maker106}) are the same, a contradiction to $\alpha=\bar{\alpha}$. \\
{\bf Case 2.3.} $D^{k}(u)=\alpha^{2}$. Then by (\ref{eq:maker68}),
we obtain

\quad $(\alpha^{2} x)^{3}\alpha^{2}=D^{k}(u^{4})\in
C^{\infty}_{1,4}$, where $x\in D_{1,4}$.\\
A direct verification leads to
\begin{eqnarray}
D^{k}(u)&=&11;\label{eq:maker107}\\
D^{k}(u^{4})&=&(1^{4}4)^{3}1^{2};\label{eq:maker108}\\
&=&1^{2}(41^{4})^{3};\label{eq:maker109}\\
&=&1^{3}(41^{4})^{2}41^{3};\label{eq:maker110}
\end{eqnarray}
or
\begin{eqnarray}
D^{k}(u)&=&44;\label{eq:maker111}\\
D^{k}(u^{4})&=&4^{3}(14^{4})^{2}14^{3}.\label{eq:maker112}
\end{eqnarray}
{\bf Case 2.3.1.} If (\ref{eq:maker107}) and (\ref{eq:maker108}), or
(\ref{eq:maker109}), or (\ref{eq:maker110}) hold, then we obtain
\begin{eqnarray}
D^{k-1}(u)&=&\alpha^{i}\bar{\alpha}\alpha\bar{\alpha}^{j},\,1\leq i,\,j\leq 3;\label{eq:maker113}\\
D^{k-1}(u^{4})&=&\alpha^{i}\bar{\alpha}\alpha\bar{\alpha}\alpha\bar{\alpha}^{4}
\alpha\bar{\alpha}\alpha\bar{\alpha}\alpha^{4}\bar{\alpha}\alpha\bar{\alpha}\alpha\bar{\alpha}^{4}
\alpha\bar{\alpha}\alpha^{j},\,1\leq i,\,j\leq 3;\label{eq:maker114}\\
&=&\alpha^{i}\bar{\alpha}\alpha\bar{\alpha}^{4}
\alpha\bar{\alpha}\alpha\bar{\alpha}\alpha^{4}\bar{\alpha}\alpha\bar{\alpha}\alpha\bar{\alpha}^{4}
\alpha\bar{\alpha}\alpha\bar{\alpha}\alpha^{j},\,1\leq i,\,j\leq 3;\label{eq:maker115}\\
&=&\alpha^{i}\bar{\alpha}\alpha\bar{\alpha}\alpha^{4}\bar{\alpha}
\alpha\bar{\alpha}\alpha\bar{\alpha}^{4}\alpha\bar{\alpha}\alpha\bar{\alpha}\alpha^{4}\bar{\alpha}
\alpha\bar{\alpha}\alpha^{j},\,1\leq i,\,j\leq 3.\label{eq:maker116}
\end{eqnarray}
Since  $D^{k-1}(u)$ and $D^{k-1}(u^{4})$ satisfy (\ref{eq:maker71}),
comparison of the last runs of (\ref{eq:maker113}) and
(\ref{eq:maker114}) to (\ref{eq:maker116}) gives  a contradiction to $\alpha^j=\bar{\alpha}^j$.\\
{\bf Case 2.3.2.} If (\ref{eq:maker111}) and (\ref{eq:maker112})
hold, then we obtain
\begin{eqnarray}
D^{k-1}(u)&=&\alpha^{i}\bar{\alpha}^{4}\alpha^{4}\bar{\alpha}^{j},\,0\leq i,\,j\leq 3;\label{eq:maker117}\\
D^{k-1}(u^{4})&=&\alpha^{i}\bar{\alpha}^{4}\alpha^{4}\bar{\alpha}^{4}\alpha\bar{\alpha}^{4}
\alpha^{4}\bar{\alpha}^{4}\alpha^{4}\bar{\alpha}\alpha^{4}\bar{\alpha}^{4}\alpha^{4}\bar{\alpha}^{4}\alpha\bar{\alpha}^{4}
\alpha^{4}\bar{\alpha}^{4}\alpha^{j},\,0\leq i,\,j\leq
3.\label{eq:maker118}
\end{eqnarray}
An argument similar to Case 2.3.1 leads to a contradiction.\\
{\bf Case 2.4.} $D^{k}(u)=\alpha^{3}$. Then by (\ref{eq:maker68}) we
have

\quad $(\alpha^{3} x)^{3}\alpha^{3}=D^{k}(u^{4})\in
C^{\infty}_{1,4}$, where $x\in D_{1,4},$\\
which means that
\begin{eqnarray}
D^{k}(u^{4})&=&(\alpha^{4}\bar{\alpha})^{3}\alpha^{3};\label{eq:maker119}\\
&=&\alpha^{3}(\bar{\alpha}\alpha^{4})^{3}.\label{eq:maker120}
\end{eqnarray}
{\bf Case 2.4.1.} $\alpha=1$. Then since  $D^{k}(u)=1^{3}$, from
(\ref{eq:maker119}) and (\ref{eq:maker120}), we obtain
\begin{eqnarray}
D^{k-1}(u)&=&\beta^{i}\bar{\beta}\beta\bar{\beta}\beta^{j},\,1\leq i,\,j\leq 3;\label{eq:maker121}\\
D^{k-1}(u^{4})&=&\beta^{i}\bar{\beta}\beta\bar{\beta}\beta\bar{\beta}^{4}
\beta\bar{\beta}\beta\bar{\beta}\beta^{4}\bar{\beta}\beta\bar{\beta}\beta\bar{\beta}^{4}
\beta\bar{\beta}\beta\bar{\beta}^{j},\,1\leq i,\,j\leq
3;\label{eq:maker122}\\
&=&\beta^{i}\bar{\beta}\beta\bar{\beta}\beta^{4}\bar{\beta}
\beta\bar{\beta}\beta\bar{\beta}^{4}\beta\bar{\beta}\beta\bar{\beta}\beta^{4}\bar{\beta}
\beta\bar{\beta}\beta\bar{\beta}^{j},\,1\leq i,\,j\leq
3.\label{eq:maker123}
\end{eqnarray}
An argument similar to Case 2.3.1 arrives at a contradiction.\\
{\bf Case 2.4.2.} $\alpha=4$. Then  since $D^{k}(u)=4^{3}$, by
(\ref{eq:maker119}) and (\ref{eq:maker120}), we have
\begin{eqnarray}
D^{k-1}(u)&=&\beta^{i}\bar{\beta}^{4}\beta^{4}\bar{\beta}^{4}\beta^{j},\,0\leq i,\,j\leq 3;\label{eq:maker124}\\
D^{k-1}(u^{4})&=&\beta^{s}\bar{\beta}^{4}\beta^{4}\bar{\beta}^{4}\beta^{4}\bar{\beta}
\beta^{4}\bar{\beta}^{4}\beta^{4}\bar{\beta}^{4}\beta\bar{\beta}^{4}\beta^{4}\bar{\beta}^{4}\beta^{4}\bar{\beta}
\beta^{4}\bar{\beta}^{4}\beta^{4}\bar{\beta}^{t};\label{eq:maker125}\\
&=&\beta^{s}\bar{\beta}^{4}\beta^{4}\bar{\beta}^{4}\beta\bar{\beta}^{4}
\beta^{4}\bar{\beta}^{4}\beta^{4}\bar{\beta}\beta^{4}\bar{\beta}^{4}\beta^{4}\bar{\beta}^{4}\beta\bar{\beta}^{4}
\beta^{4}\bar{\beta}^{4}\beta^{4}\bar{\beta}^{t},\label{eq:maker126}
\end{eqnarray}
where $0\leq s,\,t\leq 3$. An argument similar to Case 2.3.1 leads to a contradiction.\\
{\bf Case 2.5.} $D^{k}(u)=\alpha^{4}$. Then by (\ref{eq:maker68}),
we have
\begin{equation}
(\alpha^{4} x)^{3}\alpha^{4}=D^{k}(u^{4})\in C^{\infty}_{1,4}$,
where $x\in D_{1,4}.\label{eq:maker127}
\end{equation}
A direct examination shows that there is no $x\in D_{1,4}$ such that
 (\ref{eq:maker127}) holds.\\
{\bf Case 3.} $\Sigma=\{1,b\}$, $b\geq 5$. Assume on the contrary
that there exists a smooth word $u$ of length more than 1 such that
$u^{h(1,b)}\in C^{\infty}_{1,b}$. Then $u$ has at least two runs.
Let $k$ be the maximal integer such that $D^{k-1}(u)$ has at least
two runs. Then by Power Derivative Formula, we obtain
\begin{equation}
(D^{j}(u)x)^{h(1,b)-1}D^{j}(u)=D^{j}(u^{h(1,b)})\in
C^{\infty}_{1,b}\, \textrm{for}\,1\leq j\leq k,\, \textrm{where}\,
x\in D_{1,b},\label{eq:maker128}
\end{equation}
and $D^{k}(u)$ has at most one run. Therefore

\quad $D^{k}(u)=\alpha^{i}$, where $i=0,1,\cdots,b$.\\
{\bf Case 3.1.} $D^{k}(u)=\varepsilon$. Then from
(\ref{eq:maker128}) we obtain
\begin{equation}
D^{k}(u^{h(1,b)})=x^{h(1,b)-1}\in C_{1,b}^{\infty},\;
\textrm{where}\; x\in D_{1,b}.\label{eq:maker129}
\end{equation}
{\bf Case 3.1.1.} $2\nmid b$. Then since $b\geq 5$, we have
$h(1,b)=\frac{b+5}{2}\leq b$ and $2(h(1,b)-1)-2=b+1$. By virtue of
 (\ref{eq:maker129}), a direct verification leads to
\begin{equation}
D^{k}(u^{h(1,b)})=\varepsilon,\,1^{h(1,b)-1},\,b^{h(1,b)-1}.\label{eq:maker130}
\end{equation}
Hence, on the one hand, in view of $r(w)\leq |D(w)|+2$ (see page
\pageref{huang}), from  (\ref{eq:maker130}) we obtain
\begin{equation}
r(D^{k-1}(u^{h(1,b)}))\leq |D^{k}(u^{h(1,b)})|+2\leq h(1,b)-1+2\leq
b+1.\label{eq:maker1301}
\end{equation}
On the other hand, from  (\ref{eq:maker128}) ($j=k-1$), it follows
that
\begin{equation}
D^{k-1}(u^{h(1,b)})=(D^{k-1}(u)x)^{h(1,b)-1}D^{k-1}(u)\in
C^{\infty}_{1,b},\;\textrm{where}\; x\in D_{1,b}.\label{eq:maker131}
\end{equation}
Since $D^{k-1}(u)$ has at least two runs, by (\ref{eq:maker131}), we
have

\quad$r(D^{k-1}(u^{h(1,b)}))\geq 2h(1,b)= b+5$,\\
contradicts  (\ref{eq:maker1301}).\\
{\bf Case 3.1.2.} $2\mid b$. Then since $b\geq 5$, we have
$h(1,b)=\frac{b+4}{2}< b$ and $2(h(1,b)-1)-2=b$. Thus from
(\ref{eq:maker129}), a direct verification leads to
\begin{equation}
D^{k}(u^{h(1,b)})=\varepsilon,\,1^{h(1,b)-1},\,b^{h(1,b)-1},\,(1b)^{h(1,b)-1},\,(b1)^{h(1,b)-1}.\label{eq:maker1302}
\end{equation}
By the proof of Case 3.1.1, we easily see that both $2\mid b$ and
$D^{k}(u^{h(1,b)})=\varepsilon,\,1^{h(1,b)-1}\text{ or
}\\b^{h(1,b)-1}$ also give rise to a contradiction. So, we only need
to verify the last two cases. Since $D^{k-1}(u)$ has at least two
runs, by $D^k(u)=\varepsilon$, we have
\begin{equation}
D^{k-1}(u)=\alpha^{i}\bar{\alpha}^{j},\; 1\leq i,j \leq
b-1.\label{E1}
\end{equation}
And from the last two cases of  (\ref{eq:maker1302}), we obtain
\begin{eqnarray}
D^{k-1}(u^{h(1,b)})&=&\beta^{s}(\bar{\beta}\beta^{b})^{h(1,b)-1}\bar{\beta}^{t},\; 1\leq s\leq b-1,\, 0\leq t\leq b-1;\label{E2}\\
                   &=&\beta^{s}(\bar{\beta}^{b}\beta)^{h(1,b)-1}\bar{\beta}^{t},\; 0\leq s\leq b-1,\, 1\leq t\leq b-1.\label{E3}
\end{eqnarray}
{\bf Case 3.1.2.1.} If (\ref{E1}) and (\ref{E2}) hold, then since
(\ref{E1}) and (\ref{E2}) satisfy  (\ref{eq:maker131}), we have
$\alpha=\beta,\,i=s,\,j=t=1$ and $x=\alpha^{b-s}\in D_{1,b}$. Thus,
from (\ref{E1}) and (\ref{E2}) it follows that
\begin{eqnarray}
D^{k-1}(u)&=&\alpha^{s}\bar{\alpha},\,1\leq s\leq b-1;\label{E4}\\
D^{k-1}(u^{h(1,b)})&=&\alpha^{s}(\bar{\alpha}\alpha^{b})^{h(1,b)-1}\bar{\alpha},\;
1\leq s\leq b-1.\label{E5}
\end{eqnarray}
Therefore, since $x=\alpha^{b-s}\in D_{1,b}$ and $1\leq s\leq b-1$,
we get
$b-s=1,2,3,4$.\\
{\bf Case 3.1.2.1.1.} $b-s=1$. Then $s=b-1$, and from (\ref{E4}) and
(\ref{E5}), we obtain
\begin{eqnarray}
D^{k-1}(u)&=&\alpha^{b-1}\bar{\alpha};\label{E6}\\
D^{k-1}(u^{h(1,b)})&=&\alpha^{b-1}(\bar{\alpha}\alpha^{b})^{h(1,b)-1}\bar{\alpha}.\label{E7}
\end{eqnarray}
From  (\ref{eq:maker128}) ($j=k-2$) we obtain
\begin{equation}
(D^{k-2}(u)x)^{h(1,b)-1}D^{k-2}(u)=D^{k-2}(u^{h(1,b)})\in
C^{\infty}_{1,b}, \textrm{where}\; x\in D_{1,b}.\label{E8}
\end{equation}

If $\alpha=1$  then since $D^{k-2}(u)$ and $D^{k-2}(u^{h(1,b)})$
satisfy  (\ref{E8}), by (\ref{E6}) and (\ref{E7}), we get
\begin{eqnarray}
D^{k-2}(u)&=&\beta^{i}\overbrace{\bar{\beta}\beta\cdots\bar{\beta}}^{b-1}\beta^{b}\bar{\beta}^{j},\;1\leq
i \leq b-1,\;0\leq j \leq b-1.\label{E9}
\end{eqnarray}

If $2\mid (h(1,b)-1)$ then
\begin{eqnarray}
D^{k-2}(u^{h(1,b)})=\beta^{i}\overbrace{\bar{\beta}\beta\cdots\bar{\beta}}^{b-1}(\beta^{b}\overbrace{
\bar{\beta}\cdots\beta}^{b}\bar{\beta}^{b}\overbrace{\beta\cdots\bar{\beta}}^{b})^{\frac{h(1,b)-1}{2}}
\beta^{b}\bar{\beta}^{j},\label{E10}
\end{eqnarray}
where $1\leq i \leq b-1,\;0\leq j \leq b-1$.

If $2\nmid (h(1,b)-1)$ then
\begin{eqnarray}
D^{k-2}(u^{h(1,b)})=\beta^{i}\overbrace{\bar{\beta}\beta\cdots\bar{\beta}}^{b-1}(\beta^{b}\overbrace{
\bar{\beta}\cdots\beta}^{b}\bar{\beta}^{b}\overbrace{\beta\cdots\bar{\beta}}^{b})^{\frac{h(1,b)-2}{2}}
\beta^{b}\overbrace{\bar{\beta}\cdots\beta}^{b}\bar{\beta}^{b}\beta^{j},\label{E11}
\end{eqnarray}
where $1\leq i \leq b-1,\;0\leq j \leq b-1$.  Thus, in virtue of
(\ref{E9}) and either (\ref{E10}) or (\ref{E11}), comparison of the
number of the factor $\beta^{b}$ of two sides of  (\ref{E8}) arrives
at a contradiction.

If $\alpha=b$  then similarly, by (\ref{E6}) and (\ref{E7}), we get
\begin{eqnarray}
D^{k-2}(u)&=&\beta^{i}\overbrace{\bar{\beta}^{b}\beta^{b}\cdots\bar{\beta}^{b}}^{b-1}\beta\bar{\beta}^{j},\;0\leq
i \leq b-1,\;1\leq j \leq b-1.\label{E12}
\end{eqnarray}

If $2\mid (h(1,b)-1)$ then
\begin{eqnarray}
D^{k-2}(u^{h(1,b)})=\beta^{i}\overbrace{\bar{\beta}^{b}\beta^{b}\cdots\bar{\beta}^{b}}^{b-1}(\beta\overbrace{
\bar{\beta}^{b}\cdots\beta^{b}}^{b}\bar{\beta}\overbrace{\beta^{b}\cdots\bar{\beta}^{b}}^{b})^{\frac{h(1,b)-1}{2}}
\beta\bar{\beta}^{j},\label{E13}
\end{eqnarray}
where $0\leq i \leq b-1,\;1\leq j \leq b-1$.

If $2\nmid (h(1,b)-1)$ then
\begin{eqnarray}
D^{k-2}(u^{h(1,b)})=\beta^{i}\overbrace{\bar{\beta}^{b}\beta^{b}\cdots\bar{\beta}^{b}}^{b-1}(\beta\overbrace{
\bar{\beta}^{b}\cdots\beta^{b}}^{b}\bar{\beta}\overbrace{\beta^{b}\cdots\bar{\beta}^{b}}^{b})^{\frac{h(1,b)-2}{2}}
\beta\overbrace{\bar{\beta}^{b}\cdots\beta^{b}}^{b}\bar{\beta}\beta^{j},\label{E14}
\end{eqnarray}
where $0\leq i \leq b-1,\;1\leq j \leq b-1$.  Thus, by virtue of
(\ref{E12}) and either (\ref{E13}) or (\ref{E14}), comparison of the
number of the disjoint factor
$\overbrace{\bar{\beta}^{b}\beta^{b}\cdots\bar{\beta}^{b}}^{b-1}\beta\bar{\beta}$
of two sides of  (\ref{E8}) gives rise to a contradiction.\\
{\bf Case 3.1.2.1.2.} $b-s=2$. Then $s=b-2$, and from (\ref{E4}) and
(\ref{E5}), we obtain
\begin{eqnarray}
D^{k-1}(u)&=&\alpha^{b-2}\bar{\alpha};\label{E15}\\
D^{k-1}(u^{h(1,b)})&=&\alpha^{b-2}(\bar{\alpha}\alpha^{b})^{h(1,b)-1}\bar{\alpha}.\label{E16}
\end{eqnarray}

If $\alpha=1$  then since $D^{k-2}(u)$ and $D^{k-2}(u^{h(1,b)})$
satisfy  (\ref{E8}), by (\ref{E15}) and (\ref{E16}) we get
\begin{eqnarray}
D^{k-2}(u)&=&\beta^{i}\overbrace{\bar{\beta}\beta\cdots\beta}^{b-2}\bar{\beta}^{b}\beta^{j},\;1\leq
i \leq b-1,\;0\leq j \leq b-1.\label{E17}
\end{eqnarray}

If $2\mid (h(1,b)-1)$ then
\begin{eqnarray}
D^{k-2}(u^{h(1,b)})=\beta^{i}\overbrace{\bar{\beta}\beta\cdots\beta}^{b-2}(\bar{\beta}^{b}\overbrace{
\beta\cdots\bar{\beta}}^{b}\beta^{b}\overbrace{\bar{\beta}\cdots\beta}^{b})^{\frac{h(1,b)-1}{2}}
\bar{\beta}^{b}\beta^{j},\label{E18}
\end{eqnarray}
where $1\leq i \leq b-1,\;0\leq j \leq b-1$.

If $2\nmid (h(1,b)-1)$ then
\begin{eqnarray}
D^{k-2}(u^{h(1,b)})=\beta^{i}\overbrace{\bar{\beta}\beta\cdots\beta}^{b-2}(\bar{\beta}^{b}\overbrace{
\beta\cdots\bar{\beta}}^{b}\beta^{b}\overbrace{\bar{\beta}\cdots\beta}^{b})^{\frac{h(1,b)-2}{2}}
\bar{\beta}^{b}\overbrace{
\beta\cdots\bar{\beta}}^{b}\beta^{b}\bar{\beta}^{j},\label{E19}
\end{eqnarray}
where $1\leq i \leq b-1,\;0\leq j \leq b-1$.  Thus by (\ref{E17})
and either (\ref{E18}) or (\ref{E19}), comparison of the number of
the factor $\bar{\beta}^{b}$ of two sides of  (\ref{E8}) leads to a
contradiction.

If $\alpha=b$  then analogously, by (\ref{E15}) and (\ref{E16}), we
have
\begin{eqnarray}
D^{k-2}(u)&=&\beta^{i}\overbrace{\bar{\beta}^{b}\beta^{b}\cdots\beta^{b}}^{b-2}\bar{\beta}\beta^{j},\;0\leq
i \leq b-1,\;1\leq j \leq b-1.\label{E20}
\end{eqnarray}

If $2\mid (h(1,b)-1)$ then
\begin{eqnarray}
D^{k-2}(u^{h(1,b)})=\beta^{i}\overbrace{\bar{\beta}^{b}\beta^{b}\cdots\beta^{b}}^{b-2}(\bar{\beta}\overbrace{
\beta^{b}\cdots\bar{\beta}^{b}}^{b}\beta\overbrace{\bar{\beta}^{b}\cdots\beta^{b}}^{b})^{\frac{h(1,b)-1}{2}}
\bar{\beta}\beta^{j},\label{E21}
\end{eqnarray}
where $0\leq i \leq b-1,\;1\leq j \leq b-1$.

If $2\nmid (h(1,b)-1)$ then
\begin{eqnarray}
D^{k-2}(u^{h(1,b)})=\beta^{i}\overbrace{\bar{\beta}^{b}\beta^{b}\cdots\beta^{b}}^{b-2}(\bar{\beta}\overbrace{
\beta^{b}\cdots\bar{\beta}^{b}}^{b}\beta\overbrace{\bar{\beta}^{b}\cdots\beta^{b}}^{b})^{\frac{h(1,b)-2}{2}}
\bar{\beta}\overbrace{\beta^{b}\cdots\bar{\beta}^{b}}^{b}\beta\bar{\beta}^{j},\label{E22}
\end{eqnarray}
where $0\leq i \leq b-1,\;1\leq j \leq b-1$.  Thus, in view of
(\ref{E20}) and either (\ref{E21}) or (\ref{E22}), comparison of the
number of the disjoint factor
$\overbrace{\bar{\beta}^{b}\beta^{b}\cdots\beta^{b}}^{b-2}\bar{\beta}\beta$
of two sides of   (\ref{E8}) reaches a contradiction.\\
{\bf Case 3.1.2.1.3.} $b-s=3$. Then $s=b-3$, $\alpha=1$, and from
(\ref{E4}) and (\ref{E5}), we obtain
\begin{eqnarray}
D^{k-1}(u)&=&1^{b-3}b;\label{E23}\\
D^{k-1}(u^{h(1,b)})&=&1^{b-3}(b1^{b})^{h(1,b)-1}b.\label{E24}
\end{eqnarray}

Since $D^{k-2}(u)$ and $D^{k-2}(u^{h(1,b)})$ satisfy  (\ref{E8}),
from (\ref{E23}) and (\ref{E24}), we get
\begin{eqnarray}
D^{k-2}(u)&=&\beta^{i}\overbrace{\bar{\beta}\beta\cdots\bar{\beta}}^{b-3}\beta^{b}\bar{\beta}^{j},\;1\leq
i \leq b-1,\;0\leq j \leq b-1.\label{E25}
\end{eqnarray}

If $2\mid (h(1,b)-1)$ then
\begin{eqnarray}
D^{k-2}(u^{h(1,b)})=\beta^{i}\overbrace{\bar{\beta}\beta\cdots\bar{\beta}}^{b-3}(\beta^{b}\overbrace{
\bar{\beta}\cdots\beta}^{b}\bar{\beta}^{b}\overbrace{\beta\bar{\beta}\cdots\bar{\beta}}^{b})^{\frac{h(1,b)-1}{2}}
\beta^{b}\bar{\beta}^{j},\label{E26}
\end{eqnarray}
where $1\leq i \leq b-1,\;0\leq j \leq b-1$.

If $2\nmid (h(1,b)-1)$ then
\begin{eqnarray}
D^{k-2}(u^{h(1,b)})=\beta^{i}\overbrace{\bar{\beta}\beta\cdots\bar{\beta}}^{b-3}(\beta^{b}\overbrace{
\bar{\beta}\cdots\beta}^{b}\bar{\beta}^{b}\overbrace{\beta\bar{\beta}\cdots\bar{\beta}}^{b})^{\frac{h(1,b)-2}{2}}
\beta^{b}\overbrace{
\bar{\beta}\cdots\beta}^{b}\bar{\beta}^{b}\beta^{j},\label{E27}
\end{eqnarray}
where $1\leq i \leq b-1,\;0\leq j \leq b-1$.  Thus, in view of
(\ref{E25}) and either (\ref{E26}) or (\ref{E27}), comparison of the
number of the factor $\beta^{b}$ of two sides of  (\ref{E8}) gives
rise to a contradiction. \\
{\bf Case 3.1.2.1.4.} $b-s=4$. Then $s=b-4$, $\alpha=1$, and from
(\ref{E4}) and (\ref{E5}), we obtain
\begin{eqnarray}
D^{k-1}(u)&=&1^{b-4}b;\label{E28}\\
D^{k-1}(u^{h(1,b)})&=&1^{b-4}(b1^{b})^{h(1,b)-1}b.\label{E29}
\end{eqnarray}

Since $D^{k-2}(u)$ and $D^{k-2}(u^{h(1,b)})$ satisfy  (\ref{E8}),
from (\ref{E28}) and (\ref{E29}), we get
\begin{eqnarray}
D^{k-2}(u)&=&\beta^{i}\overbrace{\bar{\beta}\beta\cdots\beta}^{b-4}\bar{\beta}^{b}\beta^{j},\;1\leq
i \leq b-1,\;0\leq j \leq b-1.\label{E30}
\end{eqnarray}

If $2\mid (h(1,b)-1)$ then
\begin{eqnarray}
D^{k-2}(u^{h(1,b)})=\beta^{i}\overbrace{\bar{\beta}\beta\cdots\beta}^{b-4}(\bar{\beta}^{b}\overbrace{
\beta\cdots\bar{\beta}}^{b}\beta^{b}\overbrace{\bar{\beta}\beta\cdots\beta}^{b})^{\frac{h(1,b)-1}{2}}
\bar{\beta}^{b}\beta^{j},\label{E31}
\end{eqnarray}
where $1\leq i \leq b-1,\;0\leq j \leq b-1$.

If $2\nmid (h(1,b)-1)$ then
\begin{eqnarray}
D^{k-2}(u^{h(1,b)})=\beta^{i}\overbrace{\bar{\beta}\beta\cdots\beta}^{b-4}(\bar{\beta}^{b}\overbrace{
\beta\cdots\bar{\beta}}^{b}\beta^{b}\overbrace{\bar{\beta}\beta\cdots\beta}^{b})^{\frac{h(1,b)-2}{2}}
\bar{\beta}^{b}\overbrace{\beta\cdots\bar{\beta}}^{b}\beta^{b}\bar{\beta}^{j},\label{E32}
\end{eqnarray}
where $1\leq i \leq b-1,\;0\leq j \leq b-1$.  Thus, by virtue of
(\ref{E30}) and  either (\ref{E31}) or (\ref{E32}), comparison of
the number of the factor $\bar{\beta}^{b}$ of two sides of
(\ref{E8}) arrives at a contradiction. \\
{\bf Case 3.1.2.2.} If  (\ref{E3}) holds, then since
$\widetilde{D_{1,b}}=D_{1,b}$, by (\ref{E1}) and (\ref{E3}), we
easily see that $D^{k-1}(\tilde{u})$ and
$D^{k-1}(\tilde{u}^{h(1,b)})$ are of the form (\ref{E1}), (\ref{E2})
and satisfy (\ref{eq:maker131}). Thus, by Case 3.1.2.1, a contradiction is given rise to.\\
{\bf Case 3.2.} $D^{k}(u)=\alpha^{i}$, where $1\leq i\leq b$. Then
from  (\ref{eq:maker128}) ($j=k$) it follows that

\begin{equation}
(\alpha^{i}x)^{h(1,b)-1}\alpha^{i}\in C_{1,b}^{\infty}\;,\;
\textrm{where}\; x\in D_{1,b},\;i=1,2,\ldots,b.\label{eq:maker132}
\end{equation}
{\bf Case 3.2.1.} $|x|=1$. Then by  (\ref{eq:maker132}), we have
$x=\bar{\alpha}$  and $i=1$ or $b$, which means that

$(\alpha\bar{\alpha})^{h(1,b)-1}\alpha$ or
$(\alpha^{b}\bar{\alpha})^{h(1,b)-1}\alpha^{b}\in
C_{1,b}^{\infty}$.\\
Since
$1^{2(h(1,b)-1)-1}=D((\alpha\bar{\alpha})^{h(1,b)-1}\alpha)=D((b1)^{h(1,b)-1}b)=D^{2}((\alpha^{b}
\bar{\alpha})^{h(1,b)-1}\alpha^{b})\in C_{1,b}^{\infty}$ and
$h(1,b)\geq \frac{b+4}{2}$, so in any case, we obtain $b+1\leq
2(h(1,b)-1)-1\leq b$, a contradiction.\\
{\bf Case 3.2.2.} $|x|=2$. Then by  (\ref{eq:maker132}) and $b\geq
5$, we have

\quad $x=\bar{\alpha}\alpha,\;\alpha\bar{\alpha}$ and $i=b-1$,\\
which by  (\ref{eq:maker128}) ($j=k$),  means that
\begin{eqnarray}
D^{k}(u)&=&\alpha^{b-1};\label{eq:maker133}\\
D^{k}(u^{h(1,b)})&=&\alpha^{b-1}(\bar{\alpha}\alpha^{b})^{h(1,b)-1};\label{eq:maker134}\\
                 &=&(\alpha^{b}\bar{\alpha})^{h(1,b)-1}\alpha^{b-1}.\label{eq:maker135}
\end{eqnarray}
Since
$1^{2(h(1,b)-2)}=D^{2}((\alpha^{b}\bar{\alpha})^{h(1,b)-1}\alpha^{b-1})=D^{2}(\alpha^{b-1}
(\bar{\alpha}\alpha^{b})^{h(1,b)-1})=D^{k+2}(u^{h(1,b)})\in
C_{1,b}^{\infty}$, we obtain $2(h(1,b)-2)\leq b$. Thus if $b$ is an
odd integer, then $b+1=2(h(1,b)-2)\leq b$, a contradiction.

If $b$ is an even number then $b\geq 6$. Since if
(\ref{eq:maker133}) and (\ref{eq:maker135}) hold, then
$D^{k}(\tilde{u})\text{ and }D^{k}(\tilde{u}^{h(1,b)})$ satisfy (\ref{eq:maker133}) and (\ref{eq:maker134}).
So, we only need to check the case for $D^{k}(u)\text{ and }D^{k}(u^{h(1,b)})$ satisfing (\ref{eq:maker133}) and (\ref{eq:maker134}).\\
{\bf Case 3.2.2.1.} $\alpha=1$ and $2\mid (h(1,b)-1)$. Then since
$D^{k-1}(u)$ and $D^{k-1}(u^{h(1,b)})$ satisfy (\ref{eq:maker131}),
from (\ref{eq:maker133}) and (\ref{eq:maker134}), we obtain
\begin{eqnarray}
D^{k-1}(u)&=&\beta^{i}\overbrace{\bar{\beta}\beta\cdots\bar{\beta}}^{b-1}\beta^{j};\label{eq:maker136}\\
D^{k-1}(u^{h(1,b)})&=&\beta^{i}\overbrace{\bar{\beta}\beta\cdots\bar{\beta}}^{b-1}(\beta^{b}\overbrace{\bar{\beta}\beta\cdots\beta}^{b}
\bar{\beta}^{b}\overbrace{\beta\bar{\beta}\cdots\bar{\beta}}^{b})^{\frac{h(1,b)-1}{2}}
\beta^{j}\nonumber\\
&=&\underbrace{\beta^{i}\bar{\beta}\beta\cdots\bar{\beta}\beta^j}_{D^{k-1}(u)}\underbrace{\beta^{b-i-j}}_{x}\underbrace{\beta^i\bar{\beta}\cdots\bar{\beta}\beta}_{D^{k-1}(u)}
\bar{\beta}^{b}\cdots\bar{\beta}^{b}\underbrace{\beta\bar{\beta}\cdots\bar{\beta}\beta^j}_{D^{k-1}(u)},\label{eq:maker137}
\end{eqnarray}
where $1\leq i,j\leq b-1$.

Since (\ref{eq:maker136}) and (\ref{eq:maker137}) satisfy
(\ref{eq:maker131}), from (\ref{eq:maker137}) it follows that
$i=j=1$ and $x=\beta^{b-2}=\bar{\beta}y\bar{\beta}$, a
contradiction.\\
{\bf Case 3.2.2.2.} $\alpha=1$ and $2\nmid (h(1,b)-1)$. Then
similarly, from (\ref{eq:maker135}), we obtain
\begin{eqnarray}
D^{k-1}(u^{h(1,b)})&=&\beta^{i}\overbrace{\bar{\beta}\beta\cdots\bar{\beta}}^{b-1}(\beta^{b}\overbrace{\bar{\beta}\beta\cdots\beta}^{b}
\bar{\beta}^{b}\overbrace{\beta\bar{\beta}\cdots\bar{\beta}}^{b})^{\frac{h(1,b)-2}{2}}
\beta^{b}\overbrace{\bar{\beta}\beta\cdots\beta}^{b}\bar{\beta}^{j},\label{eq:maker139}
\end{eqnarray}
where $1\leq i,j\leq b-1$. Therefore, comparison the last run of (\ref{eq:maker136}) and (\ref{eq:maker139}) gives rise to a contradiction.\\
{\bf Case 3.2.2.3.} $\alpha=b$ and $2\mid (h(1,b)-1)$. Then since
$D^{k-1}(u)$ and $D^{k-1}(u^{h(1,b)})$ satisfy (\ref{eq:maker131}),
by (\ref{eq:maker133}) and (\ref{eq:maker134}), we get
\begin{eqnarray}
D^{k-1}(u)&=&\beta^{i}\overbrace{\bar{\beta}^{b}\beta^{b}\cdots\bar{\beta}^{b}}^{b-1}\beta^{j};\label{eq:maker141}\\
D^{k-1}(u^{h(1,b)})&=&\beta^{i}\overbrace{\bar{\beta}^{b}\beta^{b}\cdots\bar{\beta}^{b}}^{b-1}(\beta\overbrace{\bar{\beta}^{b}\beta^{b}\cdots\beta^{b}}^{b}
\bar{\beta}\overbrace{\beta^{b}\bar{\beta}^{b}\cdots\bar{\beta}^{b}}^{b})^{\frac{h(1,b)-1}{2}}
\beta^{j},\label{eq:maker142}
\end{eqnarray}
where $0\leq i,j\leq b-1$.

Since (\ref{eq:maker141}) and (\ref{eq:maker142}) satisfy
(\ref{eq:maker131}). So, by comparing the right sides of
(\ref{eq:maker141}) and (\ref{eq:maker142}), we have
$x=\beta^{1-j}\bar{\beta}^b\beta^{b-i}\in D_{1,b}$, that is, the
length of $x$ is larger than 6, contradicts the
fact that the length of the words of $D_{1,b}$ is no more than 4.\\
{\bf Case 3.2.2.4.} $\alpha=b$ and $2\nmid (h(1,b)-1)$. Then
analogously, from (\ref{eq:maker134}) we obtain
\begin{eqnarray}
D^{k-1}(u^{h(1,b)})&=&\beta^{i}\overbrace{\bar{\beta}^{b}\cdots\bar{\beta}^{b}}^{b-1}(\beta\overbrace{\bar{\beta}^{b}\cdots\beta^{b}}^{b}
\bar{\beta}\overbrace{\beta^{b}\cdots\bar{\beta}^{b}}^{b})^{\frac{h(1,b)-2}{2}}
\beta\overbrace{\bar{\beta}^{b}\cdots\beta^{b}}^{b}\bar{\beta}^{j},\label{ee}
\end{eqnarray}
where $0\leq i,j\leq b-1$. Comparison of the last runs of
(\ref{eq:maker141}) and
(\ref{ee}) leads to a contradiction.\\
{\bf Case 3.2.3.} $|x|=3$. Then by  (\ref{eq:maker132}) and $b\geq
5$, we have

\quad $x=1^{2}b,\, b1^{2}$ and i=b-2,\\
which by  (\ref{eq:maker128}) ($j=k$), means that
\begin{eqnarray}
D^{k}(u)&=&1^{b-2};\label{eq:maker144}\\
D^{k}(u^{h(1,b)})&=&(1^{b}b)^{h(1,b)-1}1^{b-2};\label{eq:maker145}\\
&=&1^{b-2}(b1^{b})^{h(1,b)-1}.\label{eq:maker146}
\end{eqnarray}
Thus, by (\ref{eq:maker145}) and (\ref{eq:maker146}), we have
$1^{2(h(1,b)-2)}=D^{2}((1^{b}b)^{h(1,b)-1}1^{b-2})=D^{2}(
(1^{b-2}\\(b1^{b})^{h(1,b)-1})=D^{k+2}(u^{h(1,b)})\in
C_{1,b}^{\infty}$, which means that $2(h(1,b)-2)\leq b$. Therefore,
if $b$ is an odd integer, then by  (\ref{e5.1}),
$b+1=2(h(1,b)-2)\leq b$, a contradiction.

If $b$ is an even number, then note that the right side of (\ref{eq:maker146}) is the reversal of the right side of (\ref{eq:maker145}),
as Case 3.2.2, we only need to verify the case for $D^{k}(u)\text{ and }D^{k}(u^{h(1,b)})$ satisfing (\ref{eq:maker144}) and (\ref{eq:maker145}).\\
{\bf Case 3.2.3.1.} $2\mid (h(1,b)-1)$. Then since $D^{k-1}(u)$ and
$D^{k-1}(u^{h(1,b)})$ satisfy  (\ref{eq:maker131}), from (
\ref{eq:maker144}) and (\ref{eq:maker145}), we obtain
\begin{eqnarray}
D^{k-1}(u)&=&\beta^{i}\overbrace{\bar{\beta}\beta\cdots\beta}^{b-2}\bar{\beta}^{j};\label{eq:maker147}\\
D^{k-1}(u^{h(1,b)})&=&\beta^{i}(\overbrace{\bar{\beta}\beta\cdots\beta}^{b}\bar{\beta}^{b}
\overbrace{\beta\bar{\beta}\cdots\bar{\beta}}^{b}\beta^{b})^{\frac{h(1,b)-1}{2}}
\overbrace{\bar{\beta}\beta\cdots\beta}^{b-2}\bar{\beta}^{j},\label{eq:maker148}
\end{eqnarray}
where $1\leq i,j\leq b-1$. Since (\ref{eq:maker147}) and
(\ref{eq:maker148}) satisfy (\ref{eq:maker131}), by comparing the
$b$-th runs of the right sides of (\ref{eq:maker147}) and
(\ref{eq:maker148}), we have $j=1$ and $x=\beta\bar{\beta}^{b}y\in
D_{1,b}$, contradicts the fact which the length of the words of
$D_{1,b}$ is less than 5.\\
{\bf Case 3.2.3.2.} $2\nmid (h(1,b)-1)$. Then similarly, from
(\ref{eq:maker145}) we obtain
\begin{eqnarray}
D^{k-1}(u^{h(1,b)})&=&\beta^{i}(\overbrace{\bar{\beta}\beta\cdots\beta}^{b}\bar{\beta}^{b}
\overbrace{\beta\bar{\beta}\cdots\bar{\beta}}^{b}\beta^{b})^{\frac{h(1,b)-2}{2}}
\overbrace{\bar{\beta}\cdots\beta}^{b}\bar{\beta}^{b}\overbrace{\beta\cdots\bar{\beta}}^{b-2}\beta^{j},\label{eq:maker150}
\end{eqnarray}
where $1\leq i,j\leq b-1$. Thus, comparison of the last runs of (\ref{eq:maker147}) and (\ref{eq:maker150}) leads to a contradiction.\\
{\bf Case 3.2.4.} $|x|=4$. Then by $x\in D_{1,b}$, we see that
$x=1^{4}$. Since $b\geq 5$, it is obvious that $(\alpha^{i}
1^{4})^{b-1}\notin C^{\infty}_{1,b}$ for $i=1,2,\ldots,b$, a contradiction.\\
{\bf Case 4.} $\Sigma=\{2,b\}$. Suppose to the contrary that there
exists a smooth word of length more than 1 such that $u^{h(2,b)}\in
C^{\infty}_{2,b}$. Then $u$ has at least two runs. Let $k$ be the
maximal integer such that $D^{k-1}(u)$ has at least two runs. Then
by Power Derivative Formula, we obtain
\begin{equation}
(D^{j}(u)x)^{h(2,b)-1}D^{j}(u)=D^{j}(u^{h(2,b)})\in
C^{\infty}_{2,b}\, \textrm{for}\,1\leq j\leq k,\, \textrm{where}\,
x\in D_{2,b},\label{eq:maker152}
\end{equation}
and $D^{k}(u)$ has at most one run. Therefore

\quad $D^{k}(u)=\alpha^{i}$, where $i=0,1,\cdots,b$.\\
{\bf Case 4.1.} $D^{k}(u)=\varepsilon$. Then by (\ref{eq:maker152}),
we have
\begin{equation}
x^{h(2,b)-1}=D^{k}(u^{h(2,b)})\in C^{\infty}_{2,b},\;
\textrm{where}\; x\in D_{2,b}.\label{eq:maker153}
\end{equation}
Note that $b\geq 3$ and $h(2,b)\geq \frac{b+3}{2}\geq 3$. From
(\ref{eq:maker153}), it follows that
\begin{equation}
D^{k}(u^{h(2,b)})=\varepsilon,\,2^{h(2,b)-1},\,b^{h(2,b)-1}.\label{eq:maker154}
\end{equation}
From (\ref{eq:maker152}) ($j=k-1$) and (\ref{eq:maker154}), we
obtain
\begin{eqnarray}
(D^{k-1}(u)x)^{h(2,b)-1}D^{k-1}(u)&=&\alpha^{i}\bar{\alpha}^{j},\,0\leq
i,j\leq b-1;\label{eq:maker155}\\
&=&\alpha^{i}\overbrace{\bar{\alpha}^{2}\alpha^{2}\cdots\beta^{2}}^{h(2,b)-1}\bar{\beta}^{j},\,1\leq
i,j\leq b-1;\label{eq:maker156}\\
&=&\alpha^{i}\overbrace{\bar{\alpha}^{b}\alpha^{b}\cdots\beta^{b}}^{h(2,b)-1}\bar{\beta}^{j},\,0\leq
i,j\leq b-1.\label{eq:maker157}
\end{eqnarray}
Since $D^{k-1}(u)$ has at least two runs, the left sides of
(\ref{eq:maker155}) to (\ref{eq:maker157}) have at least $2h(2,b)$
runs, but the right sides of (\ref{eq:maker155}) to
(\ref{eq:maker157}) have at most $h(2,b)+1$ runs, which means that
$h(2,b)\leq 1$, a
contradiction to $h(2,b)\geq 3$.\\
{\bf Case 4.2.} $D^{k}(u)=\alpha^{i},\, 1\leq i\leq b$. Then by
(\ref{eq:maker152}) ($j=k$), we have
\begin{equation}
(\alpha^{i} x)^{h(2,b)-1}\alpha^{i}=D^{k}(u^{h(2,b)})\in
C^{\infty}_{2,b},\; \textrm{where}\; x\in
D_{2,b}.\label{eq:maker158}
\end{equation}
Note that $i\geq 1$ and $h(2,b)\geq 3$. From (\ref{eq:maker158}), it
follows that $x=\varepsilon$ and $i=1$, that is,
\begin{eqnarray}
D^{k}(u)&=&\alpha;\label{eq:maker159}\\
D^{k}(u^{h(2,b)})&=&\alpha^{h(2,b)}.\label{eq:maker160}
\end{eqnarray}
Thus, on the one hand, by  (\ref{eq:maker160}), we see that
$D^{k-1}(u^{h(2,b)})$ has at most $h(2,b)+2$ runs. On the other
hand, since $D^{k-1}(u)$ has at least two runs, by
(\ref{eq:maker152}) ($j=k-1$), we see that
$D^{k-1}(u^{h(2,b)})=(D^{k-1}(u)x)^{h(2,b)-1}D^{k-1}(u)$ has at
least $2h(2,b)$ runs. Therefore, $h(2,b)\leq 2$, contradicts $h(2,b)\geq 3$.\\
{\bf Case 5.} $\Sigma=\{a,b\}$, where $b>a\geq 3$. Suppose on the
contrary that there exists a smooth word of length more than 1 such
that $u^{h(a,b)}\in C^{\infty}_{a,b}$. Then $u$ has at least two
runs. Let $k$ be the maximal integer such that $D^{k-1}(u)$ has at
least two runs. Then by Power Derivative Formula, we obtain
\begin{equation}
(D^{j}(u)x)^{h(a,b)-1}D^{j}(u)=D^{j}(u^{h(a,b)})\in
C^{\infty}_{a,b}\,\textrm{for}\,1\leq j\leq k,\,\textrm{where}\,
x\in D_{a,b},\label{eq:maker161}
\end{equation}
and $D^{k}(u)$ has at most one run. Therefore

\quad $D^{k}(u)=\alpha^{i}$, where $i=0,1,\cdots,b$.\\
{\bf Case 5.1.} $D^{k}(u)=\varepsilon$. Then by (\ref{eq:maker161})
($j=k$), we have
\begin{equation}
D^{k}(u^{h(a,b)})=x^{h(a,b)-1}\in C^{\infty}_{a,b},\;
\textrm{where}\; x\in D_{a,b}.\label{eq:maker162}
\end{equation}
Note that $b\geq 4$ and $h(a,b)\geq 4$. From (\ref{eq:maker162}), it
follows that
\begin{equation}
D^{k}(u^{h(a,b)})=\varepsilon,\,a^{h(a,b)-1},\,b^{h(a,b)-1}.\label{eq:maker163}
\end{equation}
From (\ref{eq:maker161}) ($j=k-1$) and (\ref{eq:maker163}), we
obtain
\begin{eqnarray}
(D^{k-1}(u)x)^{h(a,b)-1}D^{k-1}(u)&=&\alpha^{i}\bar{\alpha}^{j},\,0\leq
i,j\leq b-1;\label{eq:maker164}\\
&=&\alpha^{i}\overbrace{\bar{\alpha}^{a}\alpha^{a}\cdots\beta^{a}}^{h(a,b)-1}\bar{\beta}^{j},\,1\leq
i,j\leq b-1;\label{eq:maker165}\\
&=&\alpha^{i}\overbrace{\bar{\alpha}^{b}\alpha^{b}\cdots\beta^{b}}^{h(a,b)-1}\bar{\beta}^{j},\,0\leq
i,j\leq b-1.\label{eq:maker167}
\end{eqnarray}
Since $D^{k-1}(u)$ has at least two runs, the left sides of
(\ref{eq:maker164}) to (\ref{eq:maker167}) have at least $2h(2,b)$
runs, but the right sides of (\ref{eq:maker164}) to
(\ref{eq:maker167}) have at most $h(a,b)+1$ runs, which means that
$h(a,b)\leq 1$, a
contradiction to $h(a,b)\geq 4$.\\
{\bf Case 5.2.} $D^{k}(u)=\alpha^{i},\, 1\leq i\leq b$. Then by
(\ref{eq:maker161}) ($j=k$), we have
\begin{equation}
(\alpha^{i} x)^{h(a,b)-1}\alpha^{i}=D^{k}(u^{h(2,b)})\in
C^{\infty}_{a,b},\,\textrm{where}\,x\in D_{a,b}.\label{eq:maker168}
\end{equation}
Note that $i\geq 1$ and $h(a,b)\geq 4$. From (\ref{eq:maker168}), it
follows that $x=\varepsilon$ and $i=1$, that is,
\begin{eqnarray}
D^{k}(u)&=&\alpha;\label{eq:maker169}\\
D^{k}(u^{h(a,b)})&=&\alpha^{h(a,b)}.\label{eq:maker170}
\end{eqnarray}
Therefore, on the one hand, by  (\ref{eq:maker170}), we see that
$D^{k-1}(u^{h(a,b)})$ has at most $h(a,b)+2$ runs. On the other
hand, since $D^{k-1}(u)$ has at least two runs, by
(\ref{eq:maker161}) ($j=k-1$), we see that
$D^{k-1}(u^{h(a,b)})=(D^{k-1}(u)x)^{h(a,b)-1}D^{k-1}(u)$ has at
least $2h(a,b)$ runs. Therefore, $h(a,b)\leq 2$, contradicts
$h(a,b)\geq 4$.

{\bf (2)} From  (\ref{e5.1}), it immediately follows that if $b\geq
5$ then $h(a,b)\leq b$, and
$h(1,2)=3,\,h(1,3)=5,\,h(1,4)=4,\,h(2,3)=3,\,h(2,4)=4,\,h(3,4)=4$.
Hence, except for the case $a=1,\,b=3$, we have $h(a,b)\leq b+1$.
Moreover, since $(1^{3}3131^{3}3)^{4}\in C^{\infty}_{1,3}$ and
$\alpha^{b}\in C^{\infty}_{a,b}$. Thus, by  (\ref{e5.2}), the
assertion (2) holds. $\Box$

\noindent{\bf Remark.} Analogously, one can give a new proof of
Proposition~\ref{p2} without machine computation.

\section{The number of smooth power words\label{s5}}
Note that if $|u|$ is an even number, then
$\mathit{first}(\Delta^{-1}_{\alpha}(u))\neq
\mathit{last}(\Delta^{-1}_{\alpha}(u))$. Thus from the definition of
smooth words, we obtain the result of the operator
$\Delta^{-1}_{\alpha}\;(\alpha=a$ or $b)$ as below.
\begin{lem}\label{p5}
Let $u\in \Sigma^{+}_{a,b}$, where $a$, $b$ have same parity, and
$u$ has even length. Then\rm

(1) \it $\Delta^{-k}_{\alpha}(u^{n})=(\Delta^{-k}_{\alpha}(u))^{n}$
and $\Delta^{-k}_{\alpha}(u)$ has even length for $\forall\,k,\,n\in
N$;\rm

(2) \it If\, $u^{n}\in C^{\infty}_{a,b}$ then
$(\Delta^{-k}_{\alpha}(u))^{n}\in C^{\infty}_{a,b}$ for
$\forall\,k,\,n\in N$.
\end{lem}
\prf{(1) Since $u$ has even length and $a,\,b$ have same parity, we
readily see that $\Delta^{-1}_{\alpha}(u)$  also has even length and
$\Delta^{-1}_{\alpha}(u^{n})=(\Delta^{-1}_{\alpha}(u))^{n}$ for
$\alpha\in \Sigma$, which suggest that the assertion (1) holds for
$k=1$.

Now suppose that the assertion (1) holds for $k=m\,(\geq 1)$. Then
we  see that $\Delta^{-m}_{\alpha}(u)$  has even length and
$\Delta^{-m}_{\alpha}(u^{n})=(\Delta^{-m}_{\alpha}(u))^{n}$, which
imply that
$\Delta^{-(m+1)}_{\alpha}(u)=\Delta^{-1}_{\alpha}(\Delta^{-m}_{\alpha}(u))$
still has even length and
\begin{eqnarray}
\Delta^{-(m+1)}_{\alpha}(u^{n})&=&\Delta^{-1}_{\alpha}(\Delta^{-m}_{\alpha}(u^{n}))\nonumber\\
                               &=&\Delta^{-1}_{\alpha}((\Delta^{-m}_{\alpha}(u))^{n})\nonumber\\
                               &=&(\Delta^{-1}_{\alpha}(\Delta^{-m}_{\alpha}(u)))^{n}\nonumber\\
                               &=&(\Delta^{-(m+1)}_{\alpha}(u))^{n}\nonumber,
\end{eqnarray}
that is, the assertion (1) also holds for $k=m+1$.

(2) Since $u^{n}=\Delta^{k}(\Delta^{-k}_{\alpha}((u^{n})))$ and
$u^{n}\in C^{\infty}_{a,b}$, by Lemma \ref{Lem4} (3), we  get
$\Delta^{-k}_{\alpha}(u^{n})\in C^{\infty}_{a,b}$}

Now we are in a position to prove the following significative
result.
\begin{thm}\label{T5}
Let $\gamma_{a,b}(n)$ denote the number of smooth words of the form
$u^{n}$ over the $2$-letter alphabet $\{a,b\}$.\rm

(1) \it If $a,\,b$ have same parity, then except for $a=1\text{ and
}b=3$,
\begin{equation}
\gamma_{a,b}(n)=\left\{\begin{array}{ll}
0 ,& n> b\\
2,&  h(a,b)\leq n\leq b\\
\infty, & n<h(a,b)
\end{array}\right.\label{e6.1};
\end{equation}\rm

(2)
\begin{equation}
\gamma_{1,3}(n)=\left\{\begin{array}{ll}
0 ,& n\geq 5\\
\infty, & n\leq 4
\end{array}\right.\label{e6.2};
\end{equation}\rm

(3) \it If $a\text{ and }b$ have the different parity, then
\begin{equation}
\gamma_{a,b}(n)=\left\{\begin{array}{ll}
0 ,& n> b\\
2, & h(a,b)\leq n\leq b\\
\infty, &  n=1
\end{array}\right.\label{e6.3};
\end{equation}\rm

(4)
\begin{equation}
\gamma_{1,2}(n)=\left\{\begin{array}{ll}
0 ,& n\geq  3\\
46,&  n=2\\
\infty, & n=1
\end{array}\right.\label{e6.1}.
\end{equation}
\end{thm}
\noindent{\bf Proof.} By Theorem \ref{T4}, we easily comprehend that
\begin{equation}
\gamma_{a,b}(n)=\left\{\begin{array}{ll}
0 ,& n\geq \delta(a,b)\\
2,&  h(a,b)\leq n< \delta(a,b)
\end{array}\right.\label{e6.5}.
\end{equation}

(1) Note that by the definition of $h(a,b)$ ( \ref{e5.1}), if $2\mid
b$ then $2(h(a,b)-1)-2=b$, if $2\nmid b$ and $a\geq 2$ then
$2(h(a,b)-1)-2=b-1$, and if $2\nmid b$ and $a=1,\,b\neq 3$ then
$2(h(a,b)-2)-1=b$. It immediately follows that if $a\geq 2$ then
$(\alpha^{a}\bar{\alpha}^{a})^{h(a,b)-1}\in C^{\infty}_{a,b}$, and
if $a=1,\,b\neq 3$ and $2\nmid b$, then $(1b^{b})^{h(a,b)-1}\in
C^{\infty}_{a,b}$. Note that if $a,\,b$ have same parity, then
$\alpha^{a}\bar{\alpha}^{a}$ and $1b^{b}$ have both even length.
Thus, by Lemma \ref{p5} (2), we get that
$(\Delta^{-k}(\alpha^{a}\bar{\alpha}^{a}))^{h(a,b)-1}$ and
$(\Delta^{-k}(1b^{b}))^{h(a,b)-1}$ are both smooth words for
$\forall\,k\in N$, which mean that $\gamma_{a,b}(n)=\infty$ for
$n<h(a,b)$. Therefore, since $\delta(a,b)=b+1$,  by (\ref{e6.5}), we
see that the assertion (1) holds.

(2) Note that $(31^{3}3131^{3})^{4}\in C^{\infty}_{1,3}$, by Lemma
\ref{p5} (2), we see that
$[\Delta^{-k}_{\alpha}(31^{3}3131^{3})]^{4}\in C^{\infty}_{1,3}$ for
every positive integer $k$, which means that
$\gamma_{1,3}(4)=\infty$.  By Theorem \ref{T4}, we have
$\delta(1,3)=h(1,3)=5$, which suggests $\gamma_{1,3}(n)=0$ for
$\forall n\geq 5$. Thus, the assertion (2) holds.

(3) If $a,b$ have different parity, then since $\delta(a,b)=b+1$,
from  (\ref{e6.5}) immediately follows the assertion (3).

(4) From Table 1 (Page 10) in Sing~\cite{Sing4} it follows that the
number of smooth square words is 46. Thus by Proposition \ref{p2} we
arrive at the desired result. $\Box$
\section{Concluding remarks\label{s7}}

By Theorem \ref{T5}, if $a$ and $b$ have differen parity, then we do
not know the value of $\gamma_{a,b}(n)$ for $2\leq n<h(a,b)$ except
for the case $a=1\text{ and }b=2$.

{\bf Open problem.} \it Compute the values of $\gamma_{a,b}(n)$ for
$2\leq n<h(a,b)$, where $a$ and $b$ have differen parity.\rm

By Theorem \ref{T4}, we see that $C^{\infty}_{a,b}$-words are
$(b+1)$-power-free except for $a=1\text{ and }b=3$. The case of
larger
$k$-letter alphabets is also challenging.\\

\end{document}